\documentclass[times,final,authoryear,ntheorem]{elsarticle}

\usepackage{framed,multirow}
\usepackage[margin=1in]{geometry}

\usepackage{amsmath,amsfonts, mathrsfs}
\usepackage[amsthm]{ntheorem}
\usepackage{amssymb}
\usepackage{latexsym}
\usepackage{float}
\usepackage{chngcntr}

\usepackage{url}
\usepackage{xcolor}
\definecolor{newcolor}{rgb}{.8,.349,.1}

\usepackage[citebordercolor=white]{hyperref}

\journal{ArXiv}

\usepackage{todonotes}
\usepackage{xcolor}

\newcommand{\mydrafttext}{}
\newcommand{\drafttext}[1]{\renewcommand{\mydrafttext}{#1}}
\newboolean{draft}  
\setboolean{draft}{true}
\ifthenelse{\boolean{draft}}
{
    \newcounter{comments}
    \drafttext{{\color{red}This document is a draft.}}
    \newcommand{\shane}[1]{\addtocounter{comments}{1}{\bf \color{red}[Shane comment \thecomments: #1]}}
    \newcommand{\josh}[1]{\addtocounter{comments}{1}{\bf \color{blue}[josh comment \thecomments: #1]}}
}
{
\newcommand{\shane}[1]{}
\newcommand{\josh}[1]{}
}

\begin{document}


\begin{frontmatter}

\title{Geometry of transit orbits in the periodically-perturbed restricted three-body problem} 

\author[1]{Joshua \snm{Fitzgerald}}
\ead{joshfitz@vt.edu}
\author[1,2]{Shane D. \snm{Ross}\corref{cor1}}
\cortext[cor1]{Corresponding author}
\ead{shaneross@vt.edu}

\address[1]{Engineering Mechanics Program, Virginia Tech, Blacksburg, VA 24061, USA}
\address[2]{Aerospace and Ocean Engineering, Virginia Tech, Blacksburg, VA 24061, USA}


\begin{abstract}
In the circular restricted three-body problem, low energy transit orbits are revealed by linearizing the governing differential equations about the collinear Lagrange points. 
This procedure fails when time-periodic perturbations are considered, such as perturbation due to the sun (i.e., the bicircular problem) or orbital eccentricity of the primaries. 
For the case of a time-periodic perturbation, the Lagrange point is replaced by a periodic orbit, 
equivalently viewed as a hyperbolic-elliptic fixed point of a symplectic map (the stroboscopic Poincar\'e map). 
Transit and non-transit orbits can be identified in the discrete map about the fixed point, 
in analogy with the geometric construction of Conley and McGehee about the index-1 saddle equilibrium point in the continuous dynamical system. 
Furthermore, though the continuous time system does not conserve the Hamiltonian energy (which is time-varying), the linearized map locally conserves a time-independent
effective Hamiltonian function.
We demonstrate that the phase space geometry of transit and non-transit orbits is preserved in going from the  unperturbed to  a periodically-perturbed situation, which carries over to the full  nonlinear equations.
	    
\end{abstract}


\end{frontmatter}


\section{Introduction}
\label{sec1}
In recent decades, investigations of the circular restricted three-body problem (CR3BP) from a dynamical systems point of view have revealed an intricate fabric of manifolds woven between planets and moons \citep{Conley1968, Conley1969, McGehee1969, LlMaSi1985,KoLoMaRo2001,JaRoLoMaFaUz2002,astakhov2004capture,GoKoLoMaMaRo2004,DeJuLoMaPaPrRoTh2005,Ross2006,RoSc2007,Reddy2008,GaMaDuCa2009,topputo2013optimal,oshima2014applications,onozaki2017tube,todorovic2020arches,ren2012numerical}. These manifolds separate low-energy \textit{transit trajectories} that successfully pass through neck regions of permitted motion about the Lagrange points, thereby travelling between phase space realms of interest, from \textit{non-transit trajectories} that fail to pass through the neck regions. The phase space structures that separate transit and non-transit trajectories appear when linearizing the governing differential equations about the system's equilibria in the co-orbiting (rotating) frame, the collinear Lagrange points (particularly $L_1$ and $L_2$). 
Linearization nonetheless fails on generalizations of the circular restricted three-body problem subject to time-dependent perturbations, such as fourth-body effects (i.e., the bicircular problem) or orbital eccentricity of the primaries, because the fixed Lagrange points are no longer equilibria. 
Moreover, the instantaneous (moving) null points of the time varying vector field are not trajectories \citep{Wiggins2003}.

\begin{figure}
\centering
\includegraphics[width=0.7\textwidth]{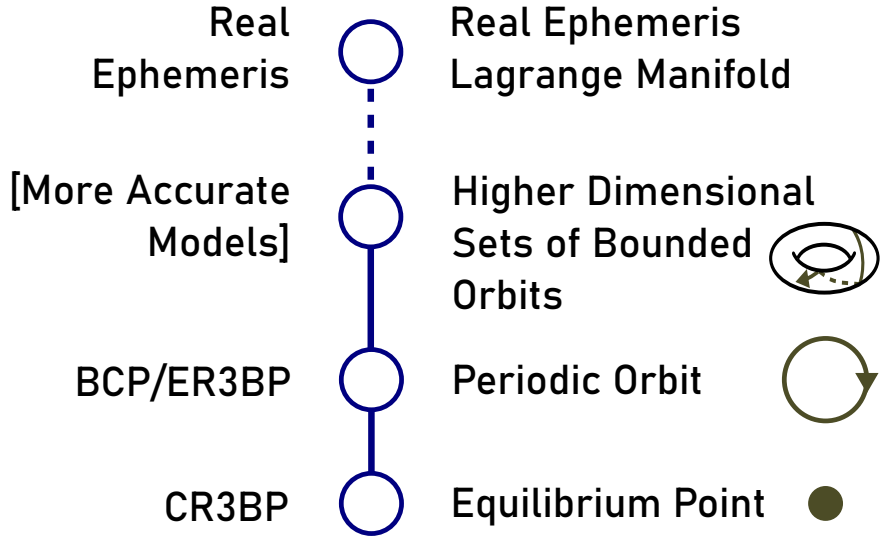}
\caption{
Schematic illustrating how the Lagrange manifold bifurcates as astrodynamical models go from simplest and least accurate at the bottom, increasing in fidelity to the real ephemeris. The bifurcation discussed in this paper is the transition from the equilibrium point  to the periodic orbit.
}
\label{fig:bifurcation}
\end{figure}

In this paper, we introduce a geometric framework for analysis of transit phenomena in time-periodic restricted three-body models like the bicircular problem (BCP) or the elliptic restricted three-body problem (ER3BP) as a natural counterpart to the time-independent circular R3BP (CR3BP). 
Higher-dimensional time-dependent manifolds, which we refer to as \textit{Lagrange manifolds}\footnote{As they are higher-dimensional analogs of the Lagrange {\it points}}, dynamically replace the $L_1$ and $L_2$ points as the fundamental objects whose stable and unstable manifolds provide the template for low energy dynamical behavior near the smaller primary. 
Under a time-periodic perturbation of period $T$, the Lagrange manifold is a manifold in the phase space diffeomorphic to $S^1$, that is, a periodic orbit with a (minimal) period equal to $T$ \citep{guckenheimer2013nonlinear}. 
Additional perturbations, not considered here, 
would further alter the topology, as depicted schematically in  Figure \ref{fig:bifurcation}.

Prior investigations into models more complicated than the CR3BP 
have successfully found periodic and quasi-periodic orbits in the vicinity of former Lagrange points by employing single shooting or multiple shooting algorithms \citep{gomez2003dynamical,jorbaetal}. 
Studies have found quasi-periodic orbits on the center manifolds of these dynamical replacements \citep{jorbaetal} and have numerically demonstrated associated transit phenomena \citep{jorba2020transport,Paez_2021}.


In this paper, we demonstrate that the linear dynamics corresponding to transit and non-transit behavior in $T$-periodically-perturbed versions of the CR3BP can be reduced to a linear time-$T$ map with the same dynamics and geometry as that in the unperturbed CR3BP.
{\it This is a significant simplification for understanding the geometry of transit orbits}, as results from several decades ago carry over in a straightforward manner, without requiring higher-order expansions.
In the phase space of the map,
the Lagrange manifold periodic orbit corresponds to an index-1  fixed point with a 1-dimensional stable manifold and 1-dimensional unstable manifold. 
Construction of transit and non-transit orbits follows from established methods dating to Conley in the 1960s \citep{Conley1968,Conley1969}.
The geometry in the linearized regime extends to the full nonlinear system, where the linear symplectic map near the Lagrange manifold will be replaced by a nonlinear symplectic map. 
Finding this nonlinear map is not our current goal, but is an objective for future research. 
According to a theorem by Moser, the linear map provides the basic geometric picture that carries over to the nonlinear case \citep{Moser1958,Moser1973}.
We demonstrate our results by considering transit orbits near the Earth-Moon $L_1$ cislunar point, the closest Lagrange point to Earth and a likely future hub for a space transportation system \citep{CoPe2001,lo2001lunar,alessi2019earth}.

The paper is organized as follows. Section \ref{classification} reviews the nature of phase space transit in the planar CR3BP. Section \ref{generaltheory} reviews an assortment of mathematical preliminaries, such as flow maps and state transition matrices, necessary for understanding the rest of the analysis. It also introduces the general theory of periodic orbit Lagrange manifolds and outlines a method for determining their initial conditions. 
Section 4 analyzes the local dynamics near an index-1 saddle-type fixed point of the Poincar\'e stroboscopic map (also called an elliptic-hyperbolic point in the discrete map context). 
Sections \ref{BCP} and \ref{ER3BP} apply these results to two examples of periodic perturbations of the CR3BP, illustrated in Figure \ref{fig:modeldiagram}: 
(i) the effect of the Sun's perturbation, known as  the {\it bicircular problem} (BCP) and (ii) the effect of the eccentricity in the Earth-Moon system,  the {\it elliptic R3BP} (ER3BP).
In putting these two distinct modifications of the R3BP on an equal footing, we seek to  emphasize the generality of our main result, the geometry of transit and non-transit orbits.

\section{Classification of orbits in the circular restricted three-body problem}\label{classification}
\subsection{Equations of motion}

\begin{figure}[h!]
\centering
\includegraphics[width=0.7 \textwidth]{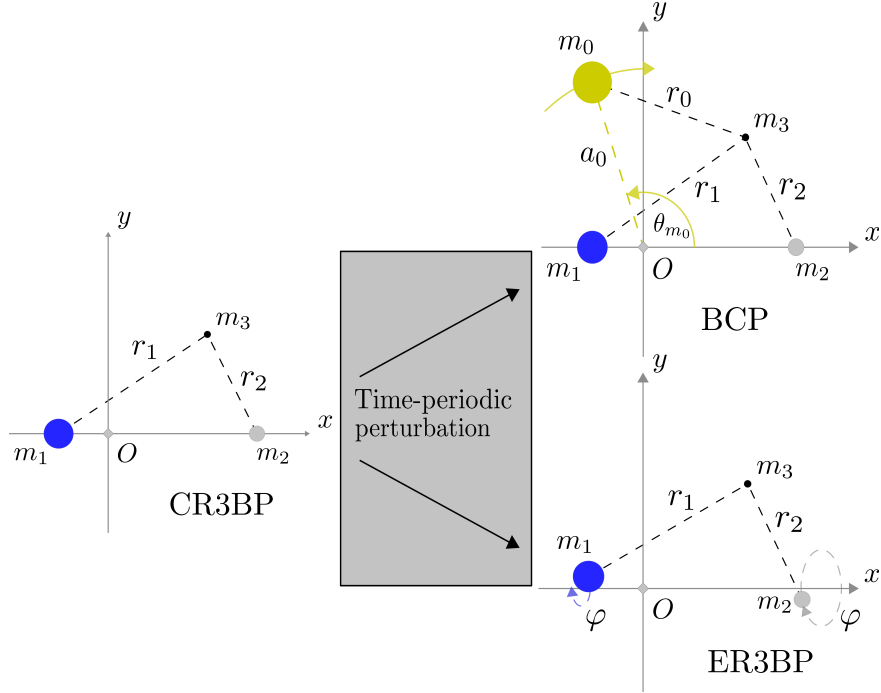}
\vspace{-4mm}\caption{The models  considered, viewed in the \(m_1\)-\(m_2\) barycentered average rotating  frame.}
\label{fig:modeldiagram}
\vspace{-4mm}
\end{figure}

The CR3BP models the motion of a small mass or test particle $m_3$ in the gravity field of two massive bodies $m_1>m_2$.
Masses \(m_1\) and \(m_2\) orbit their common center of mass \(O\) in circular orbits. 
We consider here only the planar CR3BP where  \(m_3\) is free to move throughout the $m_1$-$m_2$ orbital plane. Generalizing the following theory to the spatial CR3BP is very straightforward in the unperturbed case, and so we consider descriptions of the spatial unperturbed and perturbed cases to be beyond the scope of the current work. 
The equations of motion are written in a rotating reference frame with origin \(O\). The \(x\)-axis of the rotating frame coincides with the line between \(m_1\) and \(m_2\) whereas the \(y\)-axis points in the direction of motion of \(m_2\) 
(see Figure \ref{fig:modeldiagram}).

The non-dimensional equations of motion for \(m_3\) in the planar CR3BP (our focus here) are autonomous Hamilton's canonical equations with Hamiltonian function \citep{kolomaro},
\begin{equation}
    H_{\rm CR3BP} = 
    \tfrac{1}{2}(p_x^2+p_y^2)-x p_y + y p_x 
     - \frac{\mu_1}{r_1} - \frac{\mu_2}{r_2},
    \label{hamiltonian_pcr3bp}
\end{equation}
where,
\begin{equation}
r_1 = \sqrt{(x+\mu_2)^2+y^2}, \quad
r_2 = \sqrt{(x-\mu_1)^2+y^2},
\label{r1_r2}
\end{equation}
with $\mu_1 = 1 - \mu$ and $\mu_2 = \mu$ the non-dimensional masses of $m_1$ and $m_2$, where
$\mu = m_2/(m_1+m_2)$
is the mass parameter. 

\subsection{The Lagrange points}

The CR3BP, as an autonomous system, has five equilibrium points called \textit{Lagrange points} as viewed in the rotating frame, as shown in Figure \ref{fig:lp_pcr3bp}(a). The three equilibria lying on the \(x\)-axis, $L_1$, $L_2$, and $L_3$, are  index-1 saddle \textit{collinear points}; the remaining two, which form equilateral triangles with \(m_1\) and \(m_2\), are the 
\textit{triangular points} (center $\times$ center points for $\mu \lesssim 0.039$).
Because of their connection with low energy orbits via transit from orbits about $m_2$ and about $m_1$ and vice-versa, we focus on the collinear points.

\subsection{The Hill's region and the Hamiltonian energy}

\begin{figure}[!t]
\centering
\begin{tabular}{cc}
\includegraphics[width=0.4\textwidth]{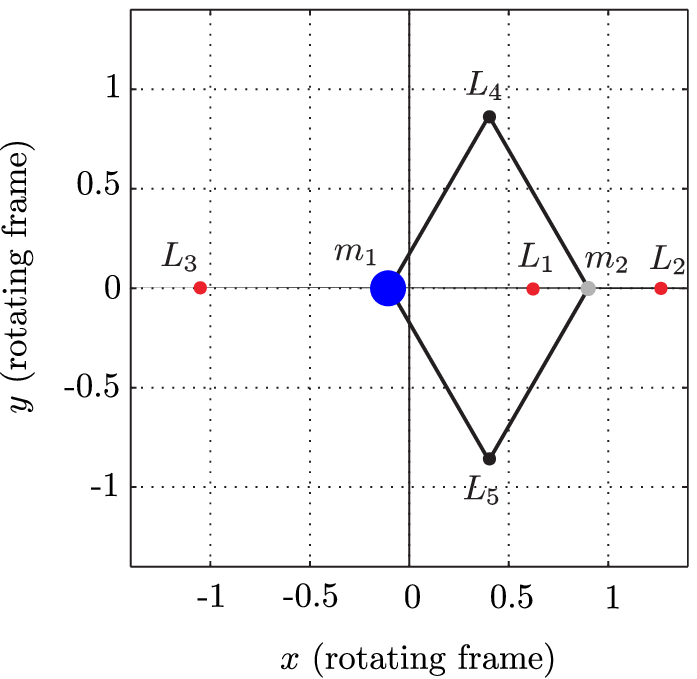} &
\includegraphics[width=0.51\textwidth]{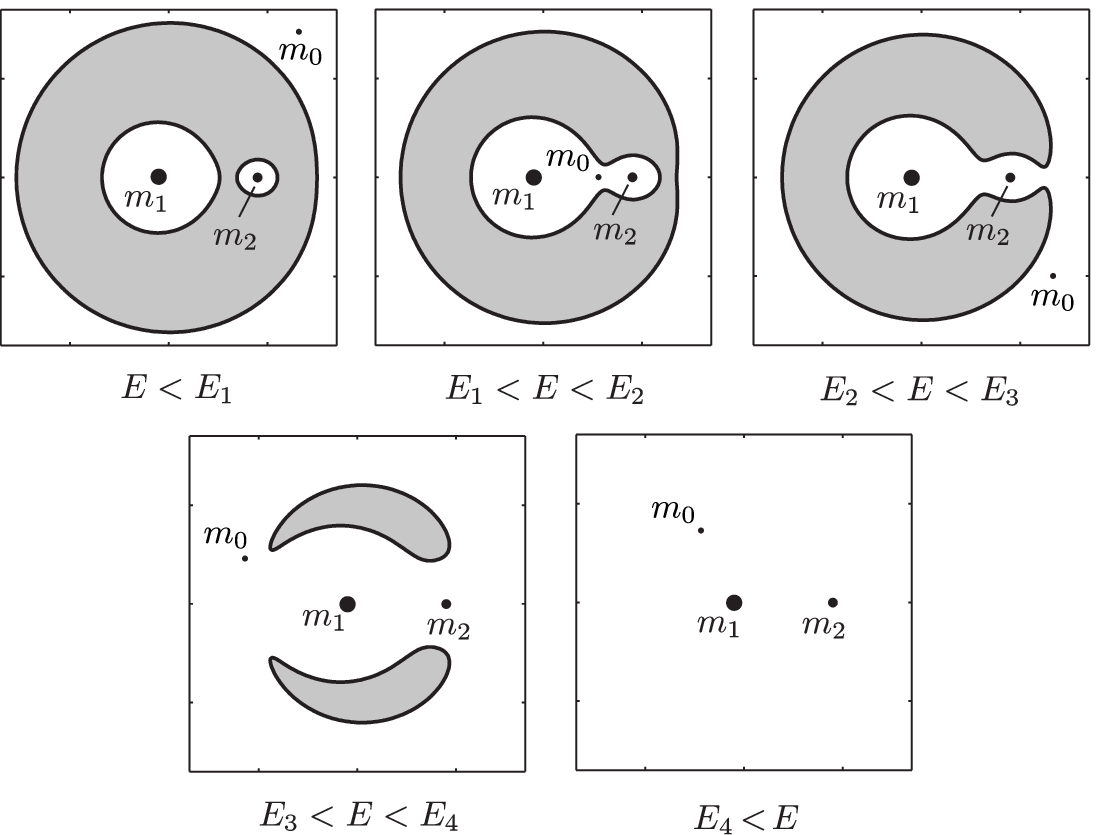}\\
{\bf (a)} &
{\bf (b)}
\end{tabular}
\caption{(a) The Lagrange points of the CR3BP for $\mu = 0.3$.
(b) The five cases of the energetically accessible regions (i.e., Hill's region) by CR3BP Hamiltonian energy. 
}
\label{fig:lp_pcr3bp}
\vspace{-4mm}
\end{figure}

Trajectories of the CR3BP conserve the Hamiltonian energy, $ H_{\rm CR3BP}=E$, where $E\in\mathbb{R}$ is the initial Hamiltonian energy.
The \textit{Hill's region} is the subset of position space throughout which $m_3$ has enough energy to travel. The boundary of the \textit{Hill's region}, beyond which lies the \textit{forbidden realm}, is called the  \textit{zero-velocity surface} in the spatial case and \textit{zero-velocity curve} in the planar case \citep{Szebehely1967}.
The qualitative characteristics of the corresponding Hill's region naturally assign $E$ to one of five different intervals 
 (see Figure \ref{fig:lp_pcr3bp}(b)):
\begin{enumerate}
    \item For $E < E_1$, $m_3$ is confined to either a subset of position space around $m_1$ (the \textit{$m_1$ realm}), a subset of position space around $m_2$ (the \textit{$m_2$ realm}), or a subset of position space outside $m_1$ and $m_2$ (the exterior realm). In this situation, $m_3$ cannot cross between any of the three realms.
    
    \item For $E_1 < E < E_2$, a {\it neck region} opens up around the $L_1$ point that permits travel between the $m_1$ and $m_2$ realms.
    
    \item For $E_2 < E < E_3$, another neck region opens up around the $L_2$ point that permits travel between the $m_2$ and exterior realms.
    
    \item For $E_3 < E < E_4$, yet another neck region opens up around the $L_3$ point that permits travel between the $m_1$ and exterior realms.
    
    \item For $E_4 < E$, the forbidden realm completely disappears.
\end{enumerate}
Thus, regions around the collinear Lagrange points play an important role in controlling transit between realms. We typically consider the second or third cases, in which transit between realms is possible but is governed by manifold structures associated with $L_1$ and in the latter case $L_2$. 

\subsection{Linearization about $L_1$ and $L_2$}

Linearizing the Hamilton's equations about $L_1$ or $L_2$,  
the eigenvalues of the linear system are a purely real pair, $\pm\lambda$, and a
purely imaginary pair, $\pm i\nu$, 
where $\lambda,\nu > 0$, which makes such points index-1 saddles \citep{MaRa1999}. The corresponding generalized
eigenvectors, when properly re-scaled, provide a {\it symplectic} eigenbasis \citep{zhong2020geometry}.
In the symplectic eigenbasis with corresponding coordinates and momenta $(q_1,p_1,q_2,p_2)$,  
the linearized equations simplify to,
\begin{equation}
    \begin{split}
        \dot{q_1} &= \lambda q_1 \text{,} \quad
        \dot{p_1} = -\lambda p_1 \text{,} \\
        \dot{q_2} &= \nu p_2 \text{,} \quad
        \dot{p_2} = -\nu q_2 \text{.} \\
    \end{split}
    \label{linearized_pcr3bp_eigen}
\end{equation}
which are Hamilton's canonical equations with corresponding quadratic Hamiltonian function, 
\begin{equation}
    H_2 = \lambda q_1 p_1 + \tfrac{1}{2}\nu(q_2^2 + p_2^2).
    \label{H2_pcr3bp}
\end{equation}
As (\ref{linearized_pcr3bp_eigen}) is linear, its solution is readily obtained and must conserve the quadratic Hamiltonian function (\ref{H2_pcr3bp}). 


\subsection{Geometry of the linearized equilibrium region}\label{four_types_orbits}

\begin{figure}[h!]
\centering
\includegraphics[width=0.7\textwidth]{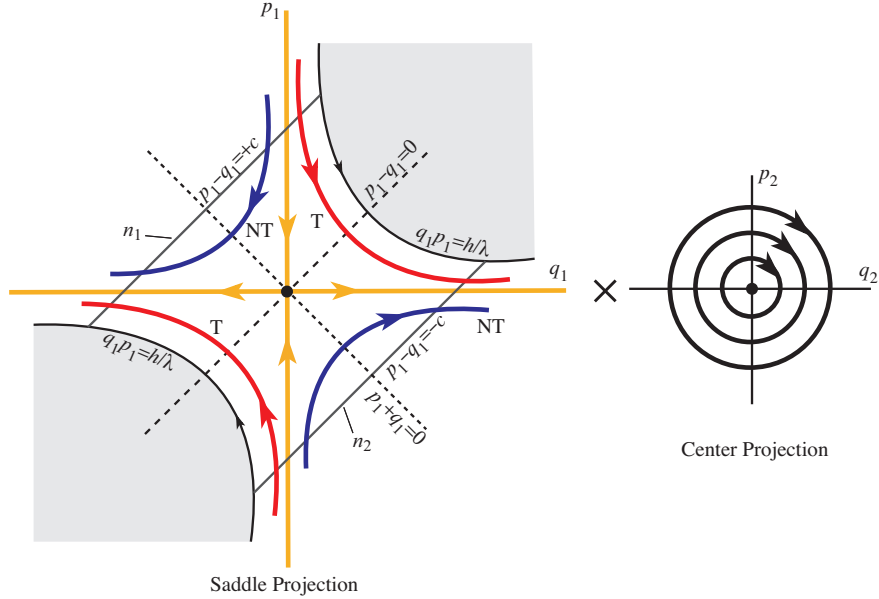}
\vspace{-3mm}
\caption{The two canonical planes of the dynamics in the symplectic eigenbasis in the neighborhood of a collinear Lagrange point; orbits labeled T transit from one realm to another, while those labeled NT do not.
}
\label{fig:linear_projections}
\end{figure}

The two canonical planes associated with (\ref{linearized_pcr3bp_eigen}) are uncoupled:\ 
the $q_1$-$p_1$ canonical plane has saddle behavior whereas the $q_2$-$p_2$ canonical plane has center behavior, as shown in
Figure \ref{fig:linear_projections}. 

Choose a fixed, small $h >0$ such that $H_2 =h$. Because $\frac{1}{2}\nu(q_2^2 + p_2^2) \ge 0$, a forbidden region in the saddle projection arises for each $h$. 
The boundary of the forbidden region is given by the hyperbolas $q_1 p_1 = h/\lambda$; 
the shape of the area outside this boundary reproduces the neck region found in the full equations of motion \citep{Conley1968}, as shown in Figure \ref{fig:linear_projections}.

For some small constant $c>0$, initial conditions along the line $p_1 - q_1 = +c$ lie entirely within one realm whereas initial conditions along the line $p_1 - q_1 = -c$ 
lie entirely within the other. For details, see \cite{kolomaro} and references therein. 
We refer to these boundaries as $n_1$ and $n_2$, respectively (see Figure \ref{fig:linear_projections}).

Orbits present in the neighborhood of the equilibrium point can be classified \citep{Conley1968} according to their behaviors in the saddle projection (see Figure \ref{fig:linear_projections}):
\begin{enumerate}
    \item The point at the origin of the saddle projection corresponds to the center manifold of the Lagrange point. Each trajectory within the center manifold is a planar periodic orbit called a \textit{Lyapunov orbit} about the equilibrium point. 
    \item The $q_1$-axis and the $p_1$-axis of the saddle projection correspond to trajectories that asymptotically approach the Lyapunov orbits as $t \to -\infty$ or $t \to + \infty$, respectively. These sets of trajectories are   the \textit{unstable and stable manifolds}, respectively, of the Lyapunov orbit of energy $h$, or, together, the \textit{asymptotic orbits}.
    \item The hyperbolic trajectories in the first and third quadrants, when integrated, intersect both $p_1 - q_1 = +c$ and $p_1 - q_1 = -c$. Because they pass from one realm to the other, they are called \textit{transit orbits}.
    \item The hyperbolic trajectories in the second and fourth quadrants are unable to intersect both $p_1 - q_1 = +c$ and $p_1 - q_1 = -c$. As they do not pass from one realm to the other, they are  
    \textit{non-transit orbits}.
\end{enumerate}
This qualitative picture in the
linearized case carries over to the nonlinear setting via 
a theorem of Moser \citep{Moser1958,Moser1973}.

\section{Lagrange manifolds in periodically-perturbed systems}\label{generaltheory}

\subsection{Periodically-perturbed systems}

In the 
analysis which follows, we consider periodically-perturbed non-autonomous dynamical systems of the form,
\begin{equation}
    \dot{\mathbf{x}} = F(\mathbf{x},t;\epsilon), \quad {\rm where} ~~\mathbf{x} \in U \subset \mathbb{R}^n, ~~t,\epsilon \in \mathbb{R}.
    \label{eq:perturbed_system}
\end{equation}
where $F$ is periodic in time $t$; that is, there exists a minimal period $T$ such that $F(\mathbf{x},t;\epsilon) = F(\mathbf{x},t + T;\epsilon)$ for all $t$, and $\epsilon$ is a perturbation parameter such that $F(\mathbf{x},t;\epsilon) \to f(\mathbf{x})$ as $\epsilon \to 0$, where $f$ is an autonomous system. A special form of $F(\mathbf{x},t;\epsilon)$ is $f(\mathbf{x}) + g(\mathbf{x},t;\epsilon)$, where $g(\mathbf{x},t;\epsilon) \to 0$ as $\epsilon \to 0$.


In a periodically-perturbed system, we can define the \textit{phase} as $\theta = \omega t \hspace{-1mm}\mod{2 \pi}$, where $\omega = 2\pi/T$. The system can then be written in 
autonomous form,
\begin{equation}
\begin{split}
    \dot{\mathbf{x}} &= F(\mathbf{x},\theta;\epsilon) \text{,} \\
    \dot{\theta} &= \omega \text{.}
    \label{perturbed_autonomous_system}
\end{split}
\end{equation}
where we note that time has been turned into a cyclic variable, $\theta \in S^1$.

\subsection{Flow maps}

Consider an arbitrary trajectory of the system (\ref{eq:perturbed_system}) with initial condition \(\mathbf{x}(t_0)=\mathbf{x}_0\). Define the corresponding flow map, $\phi(\cdot)$, as,
\begin{equation}
    \mathbf{x}(t_0) \mapsto \mathbf{x}(t)=\phi(t,t_0;\mathbf{x}_0).
\end{equation}

Consider the family of time-$T$ {\it stroboscopic} maps \(P_{t_0}:U \rightarrow U\) defined as, 
\begin{equation}
    \mathbf{x}_0  \mapsto P_{t_0}(\mathbf{x}_0) =  \phi(t_0+T,t_0;\mathbf{x}_0).
    \label{stroboscopic_map}
\end{equation} 
For a time-periodic Hamiltonian system, $P_{t_0}$ is a 
symplectic, stroboscopic map of the phase space over one period. It can equivalently be written with the parameter as the initial phase $\theta_0 = \omega t_0$ as $P_{\theta_0}$.
Note that \(P_{t_0}(\mathbf{x}_0)\) has an inverse,
\begin{equation}
\mathbf{x}_0 \mapsto P^{-1}_{t_0}(\mathbf{x}_0)  = \phi(t_0-T,t_0;\mathbf{x}_0) \text{.}
\end{equation}

\subsection{State transition and monodromy matrices}

The state transition matrix $\mathbf{\Phi}(t,t_0;\mathbf{x}_0)$ linearly approximates the flow map, 
$\phi(t,t_0;\mathbf{x}_0)$.
That is, it maps how trajectories slightly displaced from 
a reference trajectory $\overline{\mathbf{x}}(t)$ evolve 
from time $t_0$ to  $t$. 
For simplicity of notation, the dependence of the state transition matrix on its initial condition
$\mathbf{x}_0 = \overline{\mathbf{x}}(t_0)$ is suppressed.
For  
(\ref{eq:perturbed_system}), $\mathbf{\Phi}(t,t_0)$ is the solution to the initial value problem
\begin{equation}
\dot{\mathbf{\Phi}}(t,t_0) = D F(\overline{\mathbf{x}}(t),t;\epsilon) \dot{\mathbf{\Phi}}(t,t_0),
\quad \mathbf{\Phi}(t_0,t_0) = I_n,
\end{equation}
where $I_n$ is the $n\times n$ identity matrix and $DF$ is the Jacobian of $F$.

For a periodic orbit, 
the \textit{monodromy matrix} is, 
\begin{equation}\mathbf{M}_{\theta_0} = \mathbf{\Phi}(t_0+T,t_0),\end{equation}
which maps small initial displacements from the periodic orbit at phase $\theta_0$ (initial time $t_0$) to their resulting displacement after one period \citep{jordan2007nonlinear}.
For Hamiltonian systems, the monodromy matrix defines a linear symplectic map \citep{lichtenberg1992regular}. 

\subsection{Lagrange periodic orbits replace Lagrange points}\label{LagPO}

In perturbed systems where the perturbation is time-periodic and sufficiently small, equilibrium points are expected to bifurcate to periodic orbits. This result follows from the Averaging Theorem \citep{guckenheimer2013nonlinear}. 
The Lagrange points of the
CR3BP 
consequently bifurcate into periodic orbits in the presence of 
periodic perturbations.
These periodic orbits, because they dynamically replace the Lagrange points, 
by definition form a class of Lagrange manifolds. 
The behavior near a      Lagrange       point     is determined via linearization of the continuous differential equations. 
By contrast, 
the behavior near a {\it Lagrange periodic orbit} is determined via monodromy matrix calculation, which yields a discrete linear map.

A Lagrange periodic orbit  has the same period as the perturbation. 
We can 
compute a Lagrange periodic orbit by solving 
a zero-finding problem:\ choose $\mathbf{\bar x}$ that minimizes the quantity $\left|\mathbf{\bar x} - P_{0}( \mathbf{\bar x})\right|$ to within a certain tolerance (where for convenience we choose the zero phase map, $P_{0}$). 
For example, an optimization method was used to find 
the Earth-Moon $L_1$ Lagrange periodic orbit in the elliptic problem (Section 6).

To obtain periodic orbits with arbitrary perturbation sizes, we can combine this methodology with continuation. By artificially decreasing the magnitude of the perturbation to nearly zero, calculating the Lagrange manifold using the approach described, and then increasing the magnitude of the perturbation slightly and using the previous initial condition as an initial guess, it is possible to "continue" the Lagrange periodic orbit out of the Lagrange point (see Appendix \ref{appendix:c} for an example in the elliptic problem).

Unlike as in the elliptic problem, our initial condition for the bicircular problem was obtained via personal communication with the authors of \cite{jorbaetal}, who utilized a multiple-shooting and continuation method.

Example initial conditions are given in Appendix A.



         \label{step:user_drag_pt}


\section{Linear 4D symplectic map near elliptic-hyperbolic point}\label{geometry}

\subsection{Definitions}\label{definitons}

Suppose a fixed point of the time-$T$ map $P_0$ has been identified
and it is of elliptic-hyperbolic type, corresponding to a periodic orbit of saddle-center type of period $T$ of a $T$-periodic 2 degree of freedom Hamiltonian system.
Let $\mathbf{x} = (q_1,p_1,q_2,p_2)$ denote the 
displacement from the fixed point within the domain of the map $P_0$. 
The linearization of $P_0$ about the fixed point (i.e., the monodromy matrix) can be put into a symplectic eigenbasis. 
Suppose  that $(q_1,p_1,q_2,p_2)$ are coordinates with respect to this symplectic eigenbasis, where the first canonically conjugate coordinate pair $(q_1,p_1)$ corresponds to the hyperbolic (or saddle) directions and the
second canonically conjugate coordinate pair $(q_2,p_2)$ 
corresponds to the elliptic (or center) directions.
In other words, the dynamics for small $\mathbf{x}$ are given by a linear 4-dimensional symplectic map,
\begin{equation}
    \mathbf{x} \mapsto \mathbf{\Lambda} \mathbf{x}
\end{equation}
where $\mathbf{\Lambda}$ is a symplectic  matrix of the block diagonal form,
\begin{equation}
\mathbf{\Lambda} =
\begin{bmatrix}
    \sigma &   0           &  ~~0  & 0 \\
    0      &   \sigma^{-1} &  ~~0  & 0 \\
    0      &   0           &  ~~\cos \psi  & \sin \psi \\
    0      &   0           &  -\sin \psi  & \cos \psi
\end{bmatrix},
\label{eq:lambda_form}
\end{equation}
for $\sigma > 1$ and for some $\psi \in S^1$. 


\subsection{The effective quadratic Hamiltonian}

{\it Proposition 1.} 
The discrete map 
$\mathbf{x} \mapsto \mathbf{\Lambda}\mathbf{x}$
is identical to the time-$T$ map of the 
linear Hamilton's canonical equations
generated by an {\it effective} quadratic Hamiltonian,
\begin{equation}
    \Tilde{H}_2 = \Tilde{\lambda} q_1 p_1 + \tfrac{1}{2} \Tilde{\nu}( q_2^2 + p_2^2 ),
\end{equation}
where,
\begin{equation}
    \Tilde{\lambda} = \tfrac{1}{T} \ln \sigma > 0, \quad \Tilde{\nu} = \tfrac{1}{T} \psi > 0 \text{.}
    \label{effective_params}
\end{equation}
For the proof, see Appendix B. 

\subsection{Geometry of the linear map}

\begin{figure}[H]
\vspace{-0mm}
\begin{center}
\includegraphics[width=0.7\textwidth]{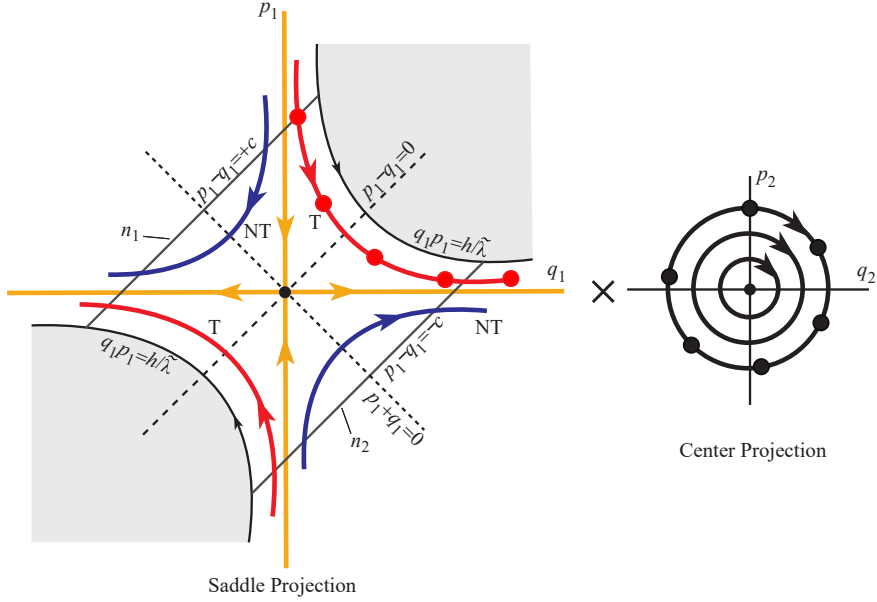}
\vspace{-3mm}
\caption{\label{fig:linear_projections_map} 
The two canonical planes of the dynamics under the mapping $\mathbf{x} \mapsto \mathbf{\Lambda} \mathbf{x}$; the orbits here are discrete  solutions of a map (represented as large dots in one of the transit curves) as compared to continuous orbits  in Figure \ref{fig:linear_projections}.}
\end{center}
\end{figure}

Because $\Tilde{H}_2$ is qualitatively identical to $H_2$ from (\ref{H2_pcr3bp}), 
the solution geometry under $\mathbf{\Lambda}$ is 
qualitatively the same as a discrete time-$T$ map of the dynamics near a collinear 
Lagrange point of the CR3BP. 
The primary difference in interpretation is  that solutions are now discrete, but still belong to families of continuous curves in the saddle and center canonical projections, 
as shown in Figure \ref{fig:linear_projections_map}.
Note that the two canonical planes are uncoupled. 
All the qualitative results related to the four types of orbits from Section \ref{four_types_orbits} carry over to the discrete case. In particular, hyperbolas in the saddle projection corresponding to transit and non-transit orbits can be identified.

\vspace{-6mm}
\subsection{Topology of the equilibrium region of the map}

\begin{figure}
\centering
\includegraphics[width=0.7 \textwidth]{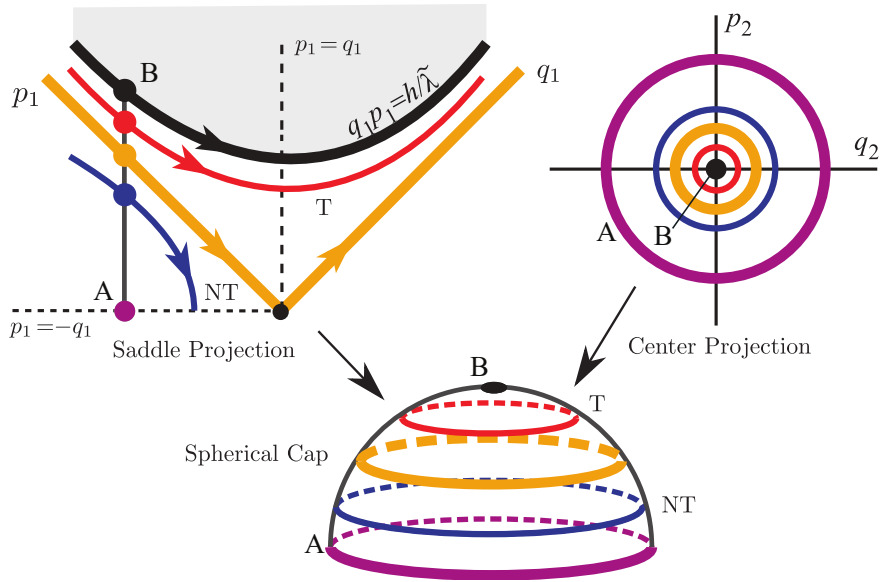}
\caption{
Construction of a  hemisphere bounding the equilibrium region of the map, along an energy manifold of energy $h$: each point along the bounding line $\overline{\rm AB}$ in the saddle projection is 
associated with a circle of initial conditions in the center projection, shrinking to a point at B.
}
\label{fig:sphericalcap}
\end{figure}

In the saddle projection, the boundaries of the equilibrium region can be defined by the two intervals of initial 
conditions parallel to the $q_1 = p_1$ line that extends between the forbidden regions. 
Pick one of the bounding intervals, say, $p_1-q_1=c$, and consider the sub-interval that enters the equilibrium region under the forward mapping, as depicted in Figure \ref{fig:sphericalcap}. 
Each point along this sub-interval corresponds to a circle in the center projection. 
The structure of the effective quadratic Hamiltonian implies that, for the trajectory on the border of the forbidden region, 
the corresponding circle shrinks to a point \citep{zhong2020geometry}. 
The bounding sub-interval is consequently homeomorphic to a spherical hemisphere; that is, $S^2 \cap \mathcal{H}^3$, where $\mathcal{H}^3$ is the upper three-dimensional half-space with boundary. 
This analysis also holds for those initial conditions that enter the region under the backward mapping, so the complete bounding interval in the saddle projection is given by $S^2$. 
Because the distance between the bounding interval and $q_1 = p_1$ is arbitrary, the entire equilibrium region is homeomorphic to $S^2 \times I$, where $I=[-c,c] \subset \mathbb{R}$ is an interval and $c>0$ is as defined in Section \ref{four_types_orbits}.

\begin{figure}
\centering
\includegraphics[width=0.35 \textwidth]{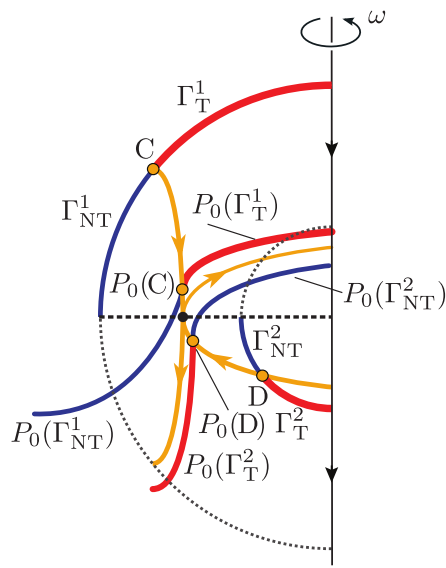}
\caption{
The McGehee representation of the discrete dynamics in the equilibrium region of the map on a fixed energy shell is obtained by rotating this diagram one revolution about the $\omega$ axis. The red lines correspond to iterates of the transit conditions under the 
stroboscopic map $P_0$; 
the blue lines correspond to iterates of the non-transit conditions; the orange lines correspond to the stable and unstable manifolds under the quadratic Hamiltonian. 
The black point corresponds to an invariant circle of the map of energy $h$, analogous to a Lyapunov orbit of energy $h$ in the unperturbed case.
}
\label{fig:mcgeheerep}
\end{figure}

The {\it McGehee representation} of the equilibrium region is 
informative for understanding 
the phase space structure of the unperturbed problem \citep{mcgehee1969some,MacKay1990,KoLoMaRo2000,waalkens2004direct,krajvnak2018phase}. 
However, we can extend the McGehee representation to the perturbed problem in a straightforward manner, as depicted in 
Figure \ref{fig:mcgeheerep}. 
The initial conditions along the boundaries of the equilibrium region that enter the region in forward time are highlighted.
For example, along the outermost bounding sphere,  $n_1$ (following the terminology of 
\cite{mcgehee1969some} and \cite{KoLoMaRo2000}), the spherical cap of transit orbits is $\Gamma_{\rm T}^1$ and the spherical band of non-transit orbits is $\Gamma_{\rm NT}^1$. 
The point C which separates $\Gamma_{\rm T}^1$ and $\Gamma_{\rm NT}^1$ is on an  orbit to an invariant circle in the equilibrium region (i.e., a quasi-periodic orbit in the full system).
The image of $\Gamma_{\rm T}^1$ and $\Gamma_{\rm NT}^1$ under the forward stroboscopic map $P_0$ is shown schematically.

Although transit initial conditions must eventually transit, they may or may not reach the other bounding sphere, $n_2$, after a single iteration of the map $P_0$, depending on their initial proximity to the stable manifold. 
Those points closest to the stable manifold will take the largest number of iterates to transit;
a discrete-time analogy to a result obtained previously (see, e.g., \cite{KoLoMaRo2000}).
Non-transit initial conditions may similarly fail to leave the equilibrium region after a single iteration. 
Thus, the transit and non-transit sets undergo stretching under the stroboscopic map.

\subsection{Connection with Lagrange periodic orbits}

A $T$-periodic Hamiltonian perturbation of the CR3BP will give rise to a Lagrange periodic orbit of period $T$ of saddle-center type.
Therefore, 
the geometry at each phase will follow the geometry given above, including in the 
{\it full nonlinear map} of the motion \citep{Wiggins2003}.

Thus, {\it the CR3BP perturbed by a periodic Hamiltonian perturbation will have the transit structure described herein.} 
Below, we consider two particular examples:\ the bicircular problem (which includes the effect of an additional mass) and the elliptic restricted three-body problem.

\section{Transit orbits in the bicircular problem}\label{BCP}

\subsection{Equations of motion of the BCP}

The bicircular problem (BCP) is a generalization of the CR3BP that describes the motions of four gravitationally interacting bodies \(m_0\), \(m_1\), \(m_2\), and \(m_3\) where \(m_2 < m_1\) and where \(m_3\) has negligible mass. 
In the inertial frame, \(m_1\) and \(m_2\) trace circular orbits around their center of mass \(O\); similarly, \(m_0\) and \(O\) trace circular orbits around their common center of mass \citep{cronin1964,SiGoJoMa1995}. 
The equations of motion are written in the CR3BP rotating reference frame so that \(m_1\) and \(m_2\) are still fixed. The large mass \(m_0\) is not fixed in the rotating frame but appears to trace a circle around \(O\) (see Figure \ref{fig:modeldiagram}).

The non-dimensional equations of motion for $m_3$ in the BCP are, unlike the equations of motion for the CR3BP, specifically time-periodic \citep{kolomaro}. 
They are Hamilton's canonical equations for a Hamiltonian,
\begin{equation}
    H_{\rm BCP} = H_{\rm CR3BP} + H_{m_0}(t),
    \label{BCP_Hamiltonian}
\end{equation}
where the time-dependent perturbation is,
\begin{equation}
    H_{m_0}(t) = \frac{\mu_0}{a_0^2} 
    \bigg(x \cos{\theta_{m_0}(t)} + y \sin{\theta_{m_0}(t)} \bigg) - \frac{\mu_0}{r_0(t)}
\end{equation}
where,
\begin{equation}
  \begin{split}
   r_0(t)^2 &= (x-a_0 \cos{\theta_{m_0}(t)})^2+(y-a_0 \sin{\theta_{m_0}(t)})^2, \\
    \theta_{m_0}(t) &= -\omega_{m_0} t + \theta_{{m_0} 0}
  \end{split}
\end{equation} 
where \(\mu_0\), \(a_0\), \(\omega_{m_0}\), \(\theta_{m_0}\), \(\theta_{{m_0} 0}\), and \(r_0\) are 
the mass, distance, angular velocity,  current angle, initial angle of \(m_0\), and distance from the particle, respectively, in non-dimensional units.
The period of \(m_0\) about the origin is $T=2\pi/\omega$ where the frequency is $\omega=\omega_{m_0}$ for this system. 
Note that the resulting equations of motion are of the form 
(\ref{eq:perturbed_system}) where $\mu_0$ corresponds to $\epsilon$.


This model has been  used to model a small celestial body or spacecraft ($m_3$) in the gravity field of the Earth ($m_1$) and Moon ($m_2$) when perturbed by the effect of the Sun ($m_0$) \citep{SiGoJoMa1995}.
The parameters in this case are \(\mu = 0.01215\), \(\mu_0 = 328900.54\), \(a_0 = 388.81114\), and \(\omega_{m_0} = 0.925195985520347\) in non-dimensional units. 

The BCP reduces to the CR3BP when gravitational perturbations from \(m_0\) are negligible; that is, when the terms due to \(m_0\) go to zero, which occurs when \(\mu_0 \to 0\) or when \(a_0 \to \infty\). The CR3BP also approximates the BCP when \(\omega_{m_0} \to \infty\) as the perturbation averages out for sufficiently large angular velocity.

\subsection{The instantaneous Lagrange points}

\begin{figure}[t!]
\begin{center}
\includegraphics[width=0.5\textwidth]{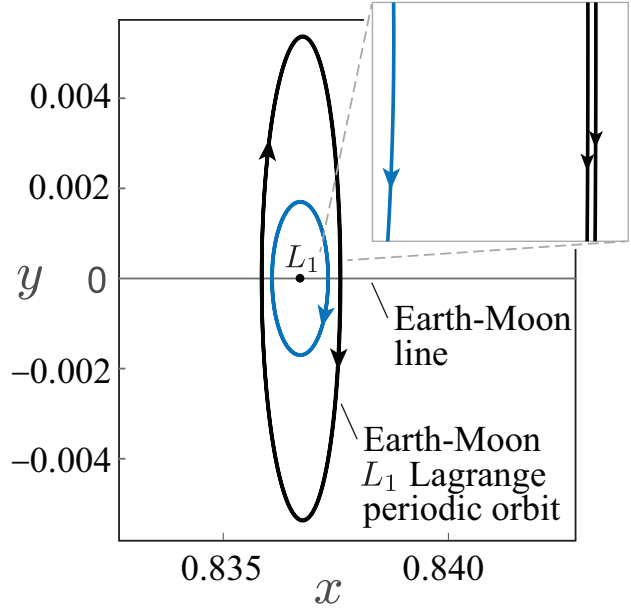}
\caption{The BCP Earth-Moon $L_1$ periodic orbit (black) compared with the path (in blue) traced by the instantaneous zero, or stagnation point, of the BCP vector field. 
The former is a trajectory; the latter is not. 
Both have doubly-looping structures over a single period of the perturbation, but at the resolution shown, even in the inset, only the periodic orbit's two loops are visible.}
\label{fig:BCML1manifold}
\end{center}
\vspace{-5mm}
\end{figure}

As discussed previously, the perturbation from \(m_0\) fundamentally removes the equilibrium points (see Figure \ref{fig:BCML1manifold}). Because the BCP is non-autonomous, the vector field associated with the equations of motion varies with \(t\) or, equivalently, \(\theta_{m_0}\). 
Setting the right side of the BCP equations of motion to zero yields an instantaneous zero of the vector field that varies 
with the independent variable, tracing out a path that repeats every  $2\pi$ in the Sun angle $\theta_{m_0}$. 
Such points are not equilibria,
and this path is \textit{not} a trajectory; particles with initial conditions along it diverge quickly.
One must consider the Lagrange periodic orbit which dynamically replaces the Lagrange point.


\subsection{Dynamics near the Sun-perturbed Earth-Moon BCP $L_1$ Lagrange periodic orbit}

The initial condition of the Sun-perturbed Earth-Moon BCP's $L_1$ Lagrange periodic orbit can be found numerically using a zero-finding procedure \citep{PaCh1989,jorbaetal}; the numerical values are given in Appendix A. 
Figure \ref{fig:BCML1manifold} depicts its path in position space.
The eigenvalues of the monodromy matrix from 0 to $T$ are
of the elliptic-hyperbolic form given previously, with
$\sigma = 4.2874\times10^{8}$ and
$\psi   = 3.0273$.
Note that the monodromy matrix could be calculated starting at a different initial phase.

The monodromy matrix of the Lagrange periodic orbit from 0 to $T$ can transformed into its symplectic eigenbasis, which is in the form of 
\eqref{eq:lambda_form}.
\begin{figure}[h!]
\centering
\begin{tabular}{cc}
\hspace{-5mm}
\includegraphics[width=0.4\textwidth]{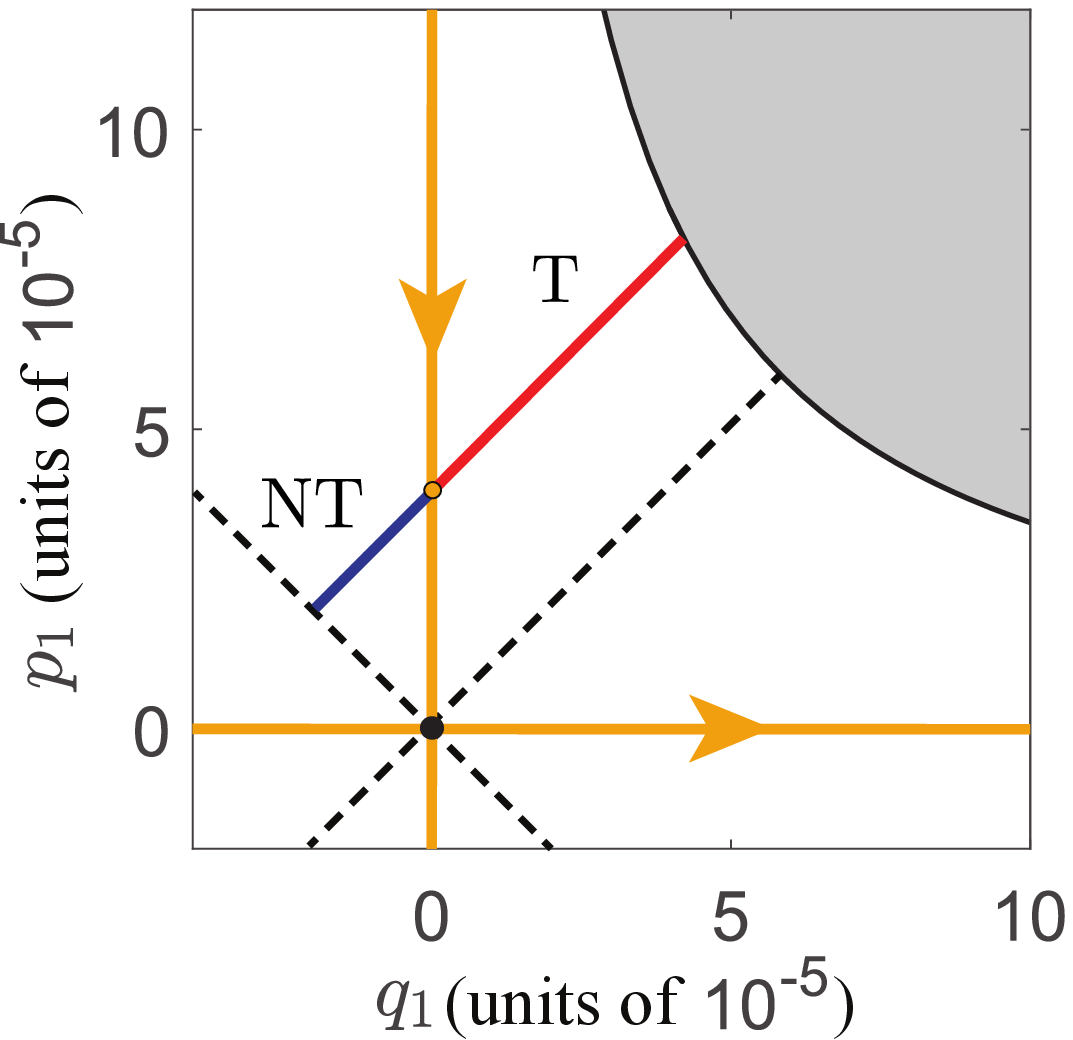} &
\includegraphics[width=0.5\textwidth]{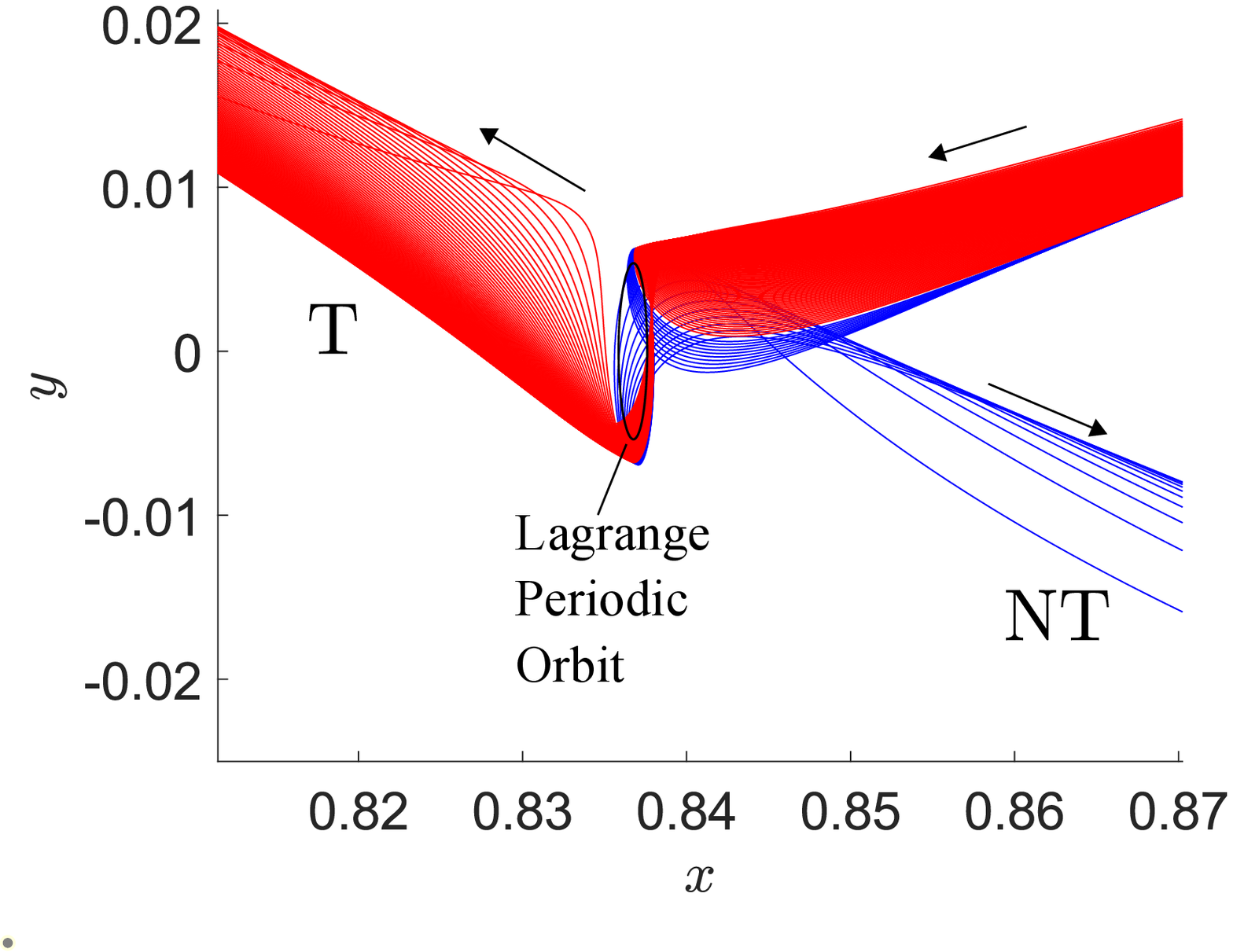}\\
{\bf (a)} &
 {\bf (b)}
\end{tabular}
\caption{(a) Numerically determined initial conditions for transit and non-transit orbits found by looking in the $q_1$-$p_1$ saddle canonical plane at initial phase $\theta=0$. $\Tilde{H}_2 = 10^{-6}$ and $c = 10^{-4}$. Compare with schematic shown in Figure \ref{fig:linear_projections_map}.
(b) The initial conditions integrated in the full equations of motion showing transit and non-transit behavior. Please refer to the online version of this article relating to color. 
}
\label{fig:dynschematicbcm}
\end{figure}
As a result, we can construct initial conditions that are transit or non-transit 
between the Earth and Moon realms when integrated in the full nonlinear equations of motion with Hamiltonian (\ref{BCP_Hamiltonian}). 
In Figure \ref{fig:dynschematicbcm}(a), the black hyperbola represents the calculated boundary of the forbidden realm, as shown schematically in Figure \ref{fig:linear_projections_map}; 
the red line contains initial conditions that should transit whereas the blue line are initial conditions that should not transit. 
In Figure \ref{fig:dynschematicbcm}(b), 
the corresponding red trajectories are transit orbits, starting in the Moon realm and going to the Earth realm, whereas the blue trajectories are non-transit orbits. 
Trajectories going from the Earth realm to Moon realm could just as easily be constructed by starting on the other boundary, $n_2$, instead of $n_1$.

\begin{figure}
\centering
\includegraphics[width=0.5 \textwidth]{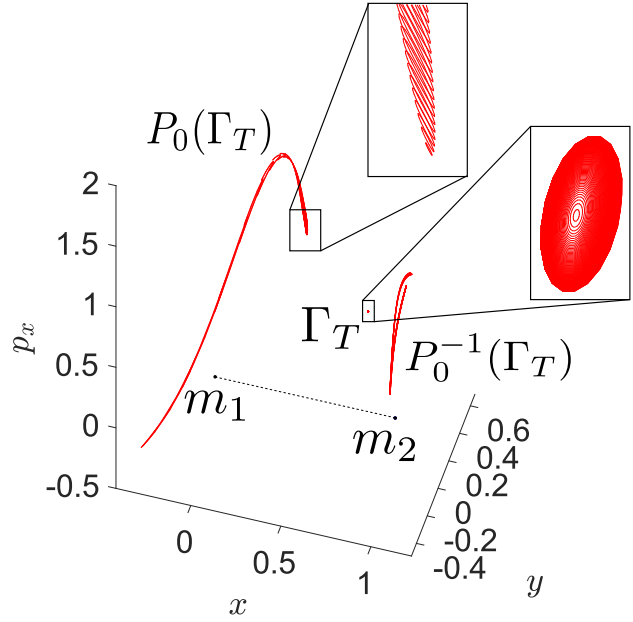}
\caption{
The spherical cap of transit orbits, $\Gamma_{\rm T}$, is mapped forwards and backwards in the bicircular model and then projected into $x$-$y$-$p_x$ space.
}
\label{fig:capmap}
\end{figure}

The spherical cap of transit orbits (labeled $\Gamma_{\rm T}$) in the bicircular model is mapped forwards and backwards for one period in Figure \ref{fig:capmap}. Under the stroboscopic map $P_0$, the set undergoes considerable distortion, but the topology, which is equivalent to that of a spherical cap, is still preserved. This setup is analogous to the description of Poincar\'e section transit orbit intersections previously computed in the Earth-Moon CR3BP \citep{KoLoMaRo2001a,deoliveira2020order}.

\section{Transit orbits in elliptic restricted three-body problem}\label{ER3BP}

\subsection{Equations of motion in the ER3BP}

The elliptic restricted three-body problem (ER3BP) is a generalization of the CR3BP that drops the restriction  that \(m_1\) and \(m_2\) move on circular orbits about their barycenter
\citep{broucke1969stability,Szebehely1967}. Instead, \(m_1\) and \(m_2\) move in more realistic elliptical orbits around their center of mass \(O\). 
We write the equations of motion in the rotating reference frame which rotates uniformly with the 
mean angular motion ($\omega =1$); 
that is, we utilize the same rotating frame as used for the CR3BP. Most authors analyzing this system utilize  a ``pulsating'' coordinate system \citep{broucke1969stability,GaMaDuCa2009}, which we have chosen not to do despite the considerable utility of this coordinate system; our aim is to bring about the commonalities of both the ER3BP and BCP and to provide ourselves with a useful toy model for our analysis.

Due to non-zero eccentricity, 
in this frame, \(m_1\) and \(m_2\) move periodically  about their CR3BP locations; their movements are given by the true anomaly $\varphi$ of the system as a function of time (see Figure \ref{fig:modeldiagram} for the geometry).
The equations of motion are Hamilton's canonical equations with Hamiltonian, 
\begin{equation}
    H_{\rm ER3BP} = \tfrac{1}{2}(p_x^2 + p_y^2) - x p_y + y p_x - \frac{\mu_1}{r_1(t)} - \frac{\mu_2}{r_2(t)},
    \label{H:er3bp}
\end{equation}
where the same non-dimensional units 
as in the CR3BP are used.
Compared to the circular problem Hamiltonian, (\ref{hamiltonian_pcr3bp}), the distances $r_i$ are now  explicit functions of time,
\begin{equation}
\begin{split}
r_i^2(t) &= 
\left| \left(\begin{bmatrix}
    x\\
    y
\end{bmatrix} + \frac{1 - \mu_i}{1+e\cos{\varphi(t)}} \mathbf{R}(t) 
\begin{bmatrix}
    \cos{\varphi(t)}\\
    \sin{\varphi(t)}
\end{bmatrix} \right) \right|^2,
\\ {\rm with}
\quad
\mathbf{R}(t) &= 
\begin{bmatrix}
   ~~\cos{t}  & \sin{t}\\
   -\sin{t} & \cos{t}
\end{bmatrix},
\end{split}
\end{equation}
where $\varphi(t)$ is the solution to the differential equation,
\begin{equation}
\dot{\varphi} = \frac{(1+e \cos{\varphi})^2}{(1-e^2)^{3/2}},
\end{equation}
with initial condition
$\varphi(0) = \varphi_0$.
For the Earth-Moon system, 
we use $e = 0.0549006$.
Using the mean anomaly as the phase $\theta$, the equations of motion are of the form 
(\ref{eq:perturbed_system}) with $T = 2 \pi/\omega=2\pi$ and with $e$ corresponding to $\epsilon$.
Note that $H_{\rm ER3BP}$ from (\ref{H:er3bp}) becomes $H_{\rm CR3BP}$ from (\ref{hamiltonian_pcr3bp}) as $e \rightarrow 0$.

\subsection{Dynamics near the Earth-Moon ER3BP $L_1$ Lagrange periodic orbit}

\begin{figure}
\begin{center}
\includegraphics[width=0.5\textwidth]{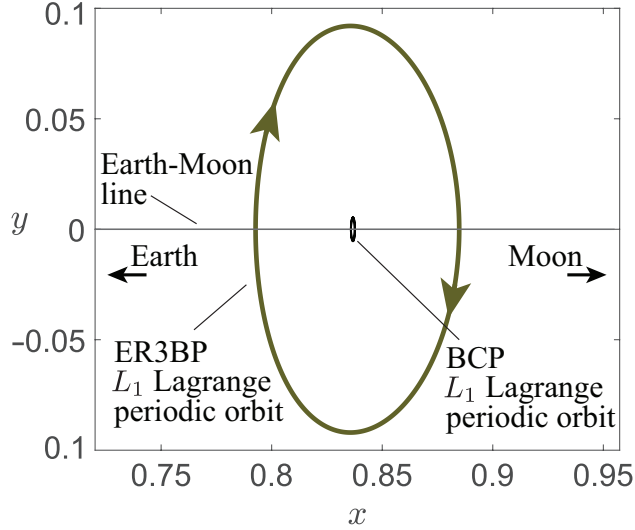}
\caption{The ER3BP Earth-Moon $L_1$ periodic orbit (large, dark green) and the BCP $L_1$ periodic orbit (black) in the position space (average rotating frame, CR3BP coordinates). The ER3BP $L_1$ periodic orbit is singly-looping, not doubly-looping as in the BCP.}
\label{fig:ER3BPL1manifold}
\end{center}
\vspace{-8mm}
\end{figure}

The initial condition of the Earth-Moon eccentric problem's $L_1$ Lagrange periodic orbit, obtained via a zero-finding algorithm (section \ref{LagPO}), is given in Appendix A. Figure \ref{fig:ER3BPL1manifold} depicts its path in position space. We show the BCP $L_1$ manifold for comparison, which is an order of magnitude smaller in amplitude.

The eigenvalues of the monodromy matrix from 0 to $T$ are
of the elliptic-hyperbolic form given in Section \ref{definitons}, with
$\sigma = 8.3659 \times 10^7$ and
$\psi   = 1.9863$.
Constructing a symplectic eigenbasis from the 
monodromy matrix 
yields initial conditions that transit or fail to transit between the Earth and Moon realms when integrated in the full nonlinear equations of 
motion---that is, Hamilton's canonical equations with Hamiltonian $H_{\rm ER3BP}$ given in \eqref{H:er3bp}. 

\begin{figure}[h!]
\centering
\begin{tabular}{cc}
\includegraphics[width=0.35\textwidth]{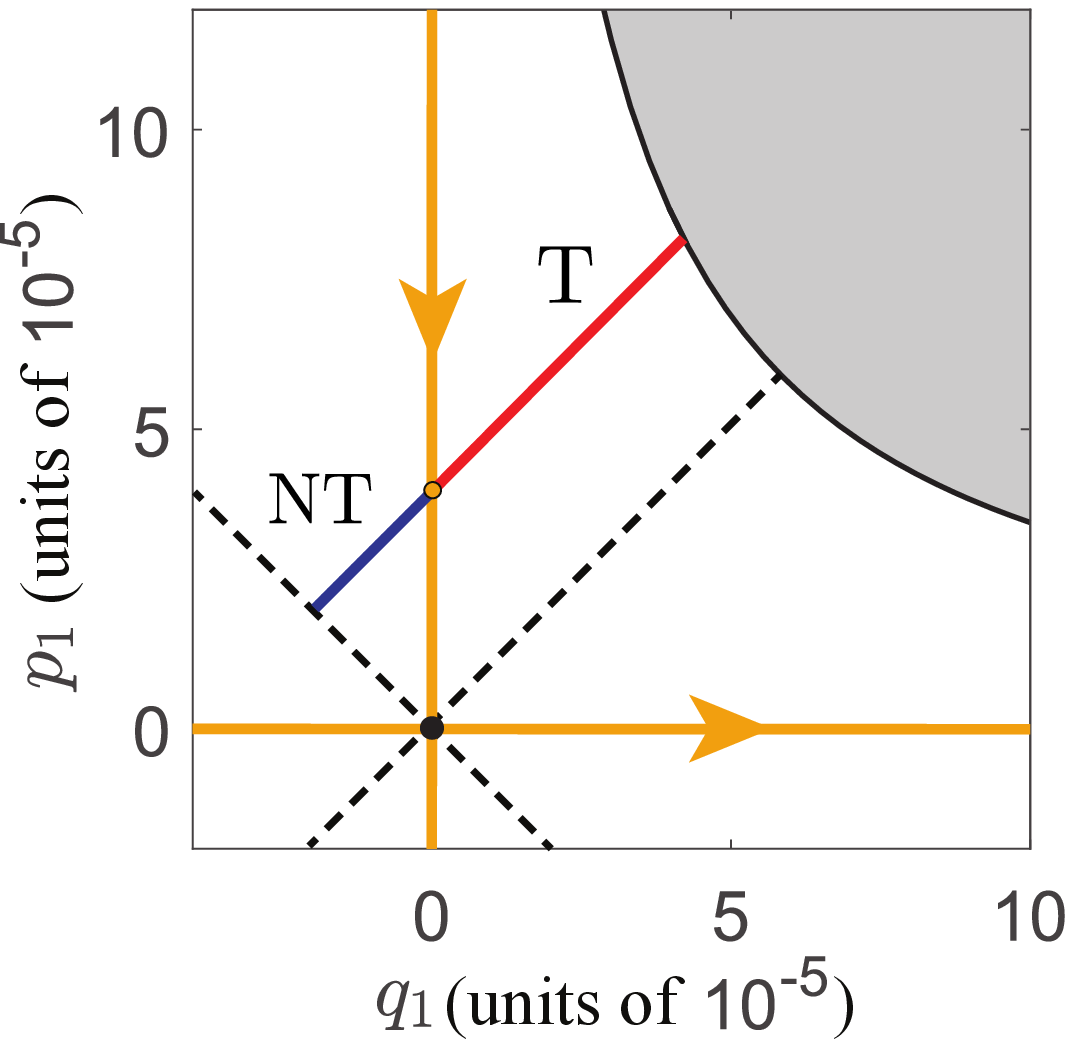} &
\hspace{-5mm}
\includegraphics[width=0.4\textwidth]{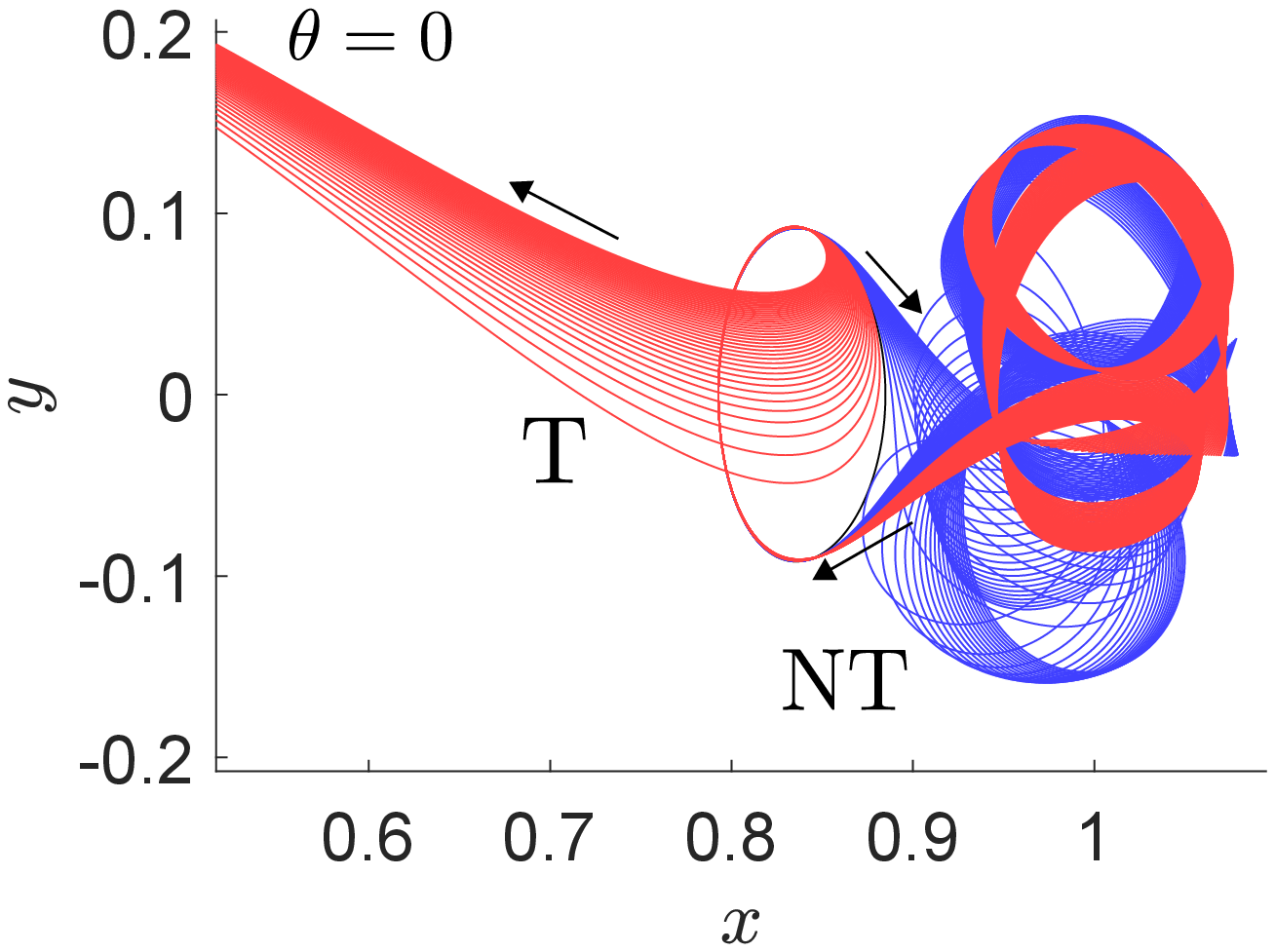}\\
{\bf (a)} & {\bf (b)}
\\
\includegraphics[width=0.4\textwidth]{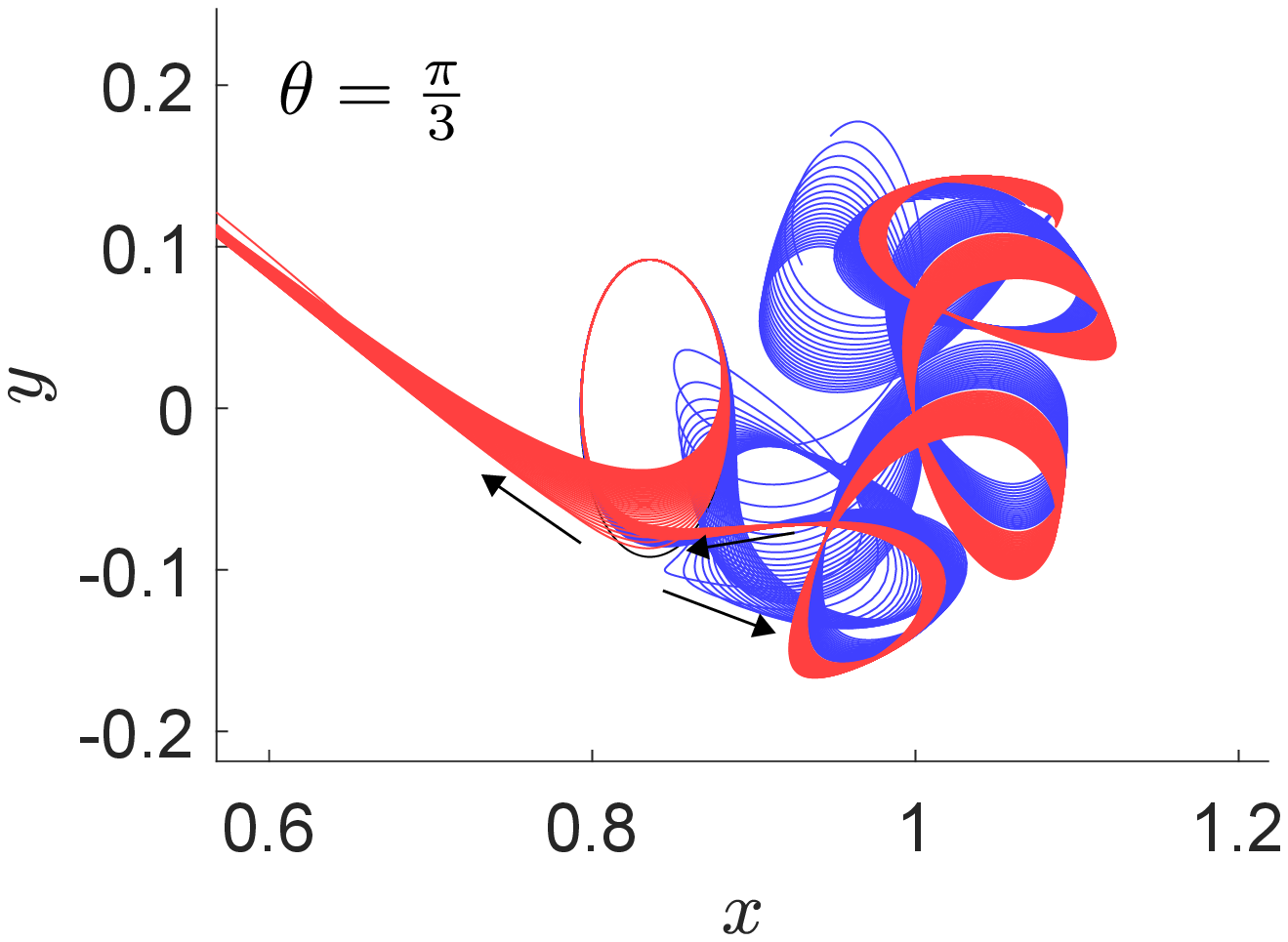} &
\hspace{-5mm}
\includegraphics[width=0.4\textwidth]{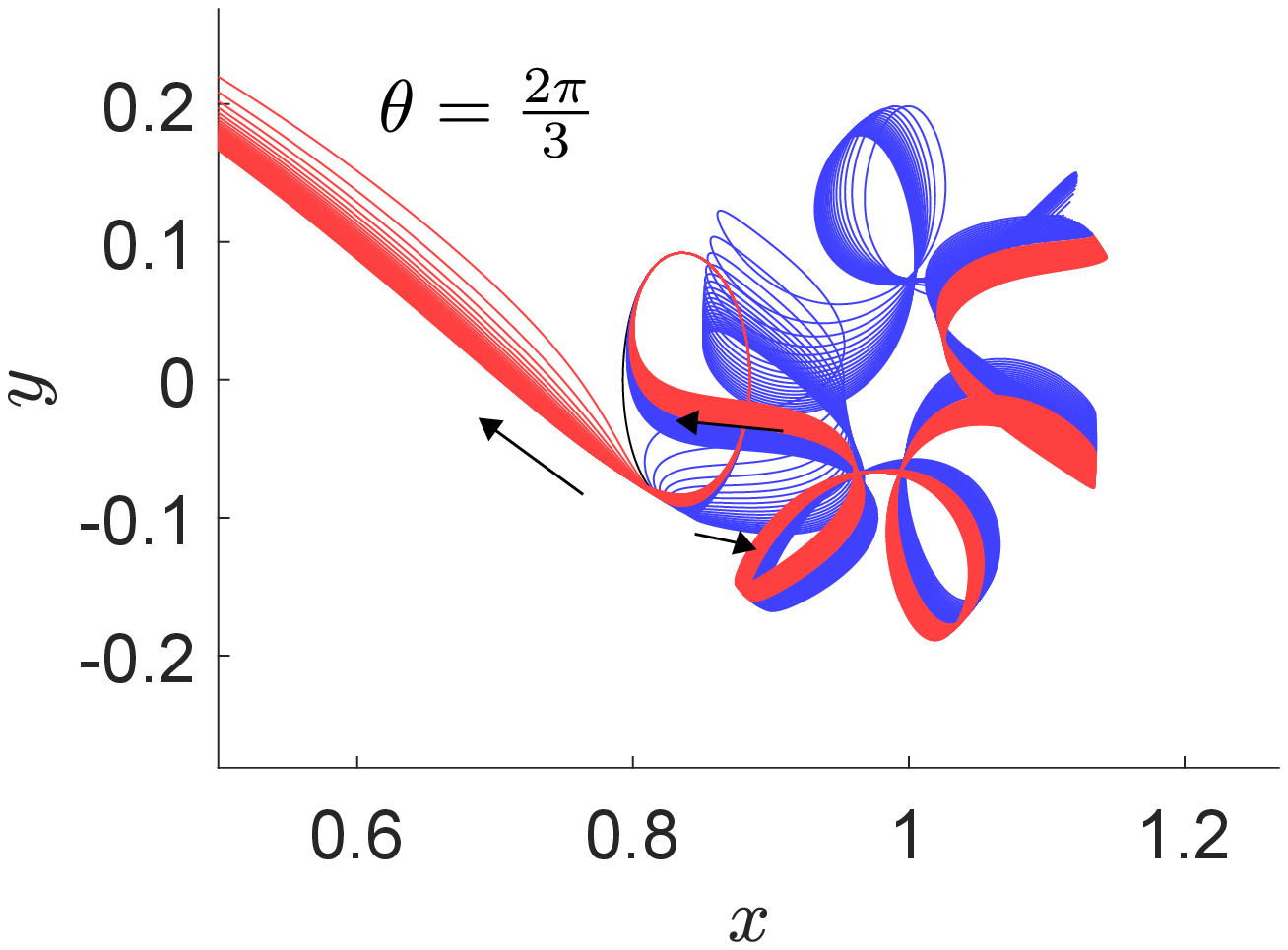}\\
{\bf (c)} & {\bf (d)}
\end{tabular}
\caption{
(a) Initial conditions for transit and non-transit orbits found by looking
in the $q_1$-$p_1$ saddle canonical plane in the symplectic eigenbasis. $\Tilde{H}_2 = 10^{-8}$ and $c = 4 \times 10^{-5}$. 
(b) The initial conditions integrated backwards and forwards in the full equations of motion, as shown, starting at 
phase (mean anomaly) $\theta=0$. 
(c) The initial conditions from part (a) integrated backwards and forwards in the full equations of motion for $\theta=\frac{\pi}{3}$. Note that the transit theory still holds at a different phase.
(d) The integrated initial conditions for $\theta=\frac{2 \pi}{3}$.  
}
\label{fig:dynschematicer3bp}
\end{figure}

In  Figure \ref{fig:dynschematicer3bp}(a), the black hyperbola represents the calculated boundary of the forbidden realm in the saddle projection. 
The red line corresponds to initial conditions, $\Gamma_{\rm T}$, that should transit whereas the blue line is initial conditions that should not transit, $\Gamma_{\rm NT}$. 
In  Figure \ref{fig:dynschematicer3bp}(b), the trajectories in the full equations of motion are shown. As expected, the red trajectories are transit orbits, starting in the Moon realm and going to the Earth realm, whereas the blue trajectories are non-transit orbits. 

Although we have shown examples of systematically finding transit and non-transit orbits for the BCP and the ER3BP at a single phase in the periodic perturbation, 
the method works equally well at other phases. We illustrate this at two additional initial phases for the initial conditions in parts (c) and (d) of Figure \ref{fig:dynschematicer3bp} for the ER3BP. 


\section{Discussion and Conclusion}

We demonstrate that the linear dynamics corresponding to transit and non-transit behavior in $T$-periodically-perturbed versions of the circular restricted three-body problem can be reduced to a linear time-$T$ map with the same orbit geometry as is now well-known in the CR3BP, going back to Conley and McGehee \citep{Conley1968,mcgehee1969some}.
Dynamically replacing the index-1 Lagrange equilibrium point of the autonomous system is a period-$T$ Lagrange periodic orbit, analyzed via a time-$T$ stroboscopic Poincar\'e map.
in the phase space of the map,
the Lagrange periodic orbit corresponds to an index-1  fixed point, or elliptic-hyperbolic point.
As we consider only the planar (two degree of freedom) problem, the Lagrange periodic orbit has a 2-dimensional center manifold, 1-dimensional stable manifold, and 1-dimensional unstable manifold.

In the extended phase space of the perturbed models (that is, including the phase of the perturbation, or cyclic time), the transit and non-transit orbits form open sets bounded by the stable and unstable manifolds to the Lagrange periodic orbit.
These results carry over to the full nonlinear system, where the linear symplectic map near the Lagrange periodic orbit is  replaced by the full nonlinear symplectic map.

Moreover, a method for elucidating the geometry of transit orbits in generalizations of the circular restricted three-body problem experiencing periodic perturbations is given. 
The Conley-McGehee representation is re-interpreted in terms of a discrete mapping rather than continuous dynamics (in Section \ref{geometry}).
The theory was demonstrated in two examples of perturbed models: the bicircular problem and the elliptic restricted three-body problem.

We illustrated our results by considering transit orbits near the Earth-Moon $L_1$ cislunar point, the most easily accessible  Lagrange point from Earth and a likely  focus for future space endeavours \citep{CoPe2001,Reddy2008,alessi2019earth,mccarthy2020cislunar,oshima2017analysis}.
Cislunar space also has significant natural connections to the Sun-Earth $L_1$ and $L_2$ regions \citep{lo2001lunar,KoLoMaRo2001a}, which can be explored using geometric techniques rather than less direct, optimization-based approaches \citep{assadian2010,onozaki2017tube,guo2019families}.

We believe that the results herein contribute significantly to the state-of-the-art in the literature. As implied in the introduction to this paper, exploring the dynamical properties of perturbations of the CR3BP has lately become a popular area of investigation in astrodynamics (refer to \citep{jorbaetal,jorba2020transport,Paez_2021,kumar2021using} for just a few recent examples). This study, by outlining a simple and straightforward method for delineating transit and nontransit behavior within perturbed models, elegantly fills an important niche in this emerging topic.

This paper also suggests a much more general discovery with ramifications beyond astrodynamics: that \textit{manifold-based transit phenomena are robust under time-periodic perturbation}. Recent scholarship has determined that manifold-based transit phemonena are also robust under dissipation \citep{zhong2020geometry, zhong2021global}. These two discoveries together help to demonstrate rigorously that natural systems subject to perturbation can exhibit the behaviors predicted by idealized tube-manifold models.

We also believe that the work herein will have useful real-world applications. Existing integrated frameworks for low-energy trajectory design utilize the dynamical characteristics of the circular restricted three-body problem \citep{kolomaro}. As shown in this paper, however, the effects of perturbations can be very large from a qualitative perspective and can permit the design of unique mission architectures that arise from the specific dynamical characteristics of perturbed models. For example, the diagrams in Sections 5 and 6 demonstrate that transit orbits "wind" on and off of the Lagrange manifolds in a way that might have practical navigational utility. 

There are several potential avenues 
for further investigation. 
This study only considered one possible topological class of Lagrange manifolds, periodic orbits generated by a single periodic perturbation. 
Additional perturbations will lead to additional bifurcations in the topology of the Lagrange point dynamical replacement
(see Figure \ref{fig:bifurcation}).
For instance, quasi-periodic Lagrange manifolds in systems with two or more perturbations of incommensurate period will generate hyperbolic structures controlling transit \citep{gomez2003dynamical,Le_Bihan_2017,jorbaetal}.

Another possibility for further study involves combining periodic perturbations with general non-conservative (e.g., dissipative, solar sail) effects \citep{bartsch2008time,zhong2020geometry}.
Our approach is applicable to the geometry of transition dynamics in other periodically-perturbed (or driven) systems governed by Hamiltonian dynamics, including chemical systems, ship dynamics,  solid state physics, and structural systems \citep{zhong2021global,NaikRoss2017a,bartsch2008time,wu2008application}.

\section*{Acknowledgments}
This work was supported in part by the National Science Foundation under awards 1537349 and 1821145.
J.F.\ was supported by a  Virginia Space Grant Consortium Graduate Research Fellowship.
The authors would like to thank Angel Jorba and Jose Rosales for providing 
an initial condition for the periodic orbit dynamical replacement to \(L_1\) in the bicircular problem.

\bibliographystyle{model5-names}
\biboptions{authoryear}
\bibliography{refs,ross_refs3}

\begin{thebibliography}{56}
\expandafter\ifx\csname natexlab\endcsname\relax\def\natexlab#1{#1}\fi
\providecommand{\url}[1]{\texttt{#1}}
\providecommand{\href}[2]{#2}
\providecommand{\path}[1]{#1}
\providecommand{\DOIprefix}{doi:}
\providecommand{\ArXivprefix}{arXiv:}
\providecommand{\URLprefix}{URL: }
\providecommand{\Pubmedprefix}{pmid:}
\providecommand{\doi}[1]{\href{http://dx.doi.org/#1}{\path{#1}}}
\providecommand{\Pubmed}[1]{\href{pmid:#1}{\path{#1}}}
\providecommand{\bibinfo}[2]{#2}
\ifx\xfnm\relax \def\xfnm[#1]{\unskip,\space#1}\fi
\bibitem[{Alessi et~al.(2019)Alessi, Masdemont \& Rossi}]{alessi2019earth}
\bibinfo{author}{Alessi, E.~M.}, \bibinfo{author}{Masdemont, J.}, \&
  \bibinfo{author}{Rossi, A.} (\bibinfo{year}{2019}).
\newblock \bibinfo{title}{The {E}arth--{M}oon system as a dynamical
  laboratory}.
\newblock {\it \bibinfo{journal}{Frontiers in Astronomy and Space Sciences}\/},
   {\it \bibinfo{volume}{6}\/}, \bibinfo{pages}{43}.
\bibitem[{Assadian \& Pourtakdoust(2010)}]{assadian2010}
\bibinfo{author}{Assadian, N.}, \& \bibinfo{author}{Pourtakdoust, S.~H.}
  (\bibinfo{year}{2010}).
\newblock \bibinfo{title}{Multiobjective genetic optimization of
  {E}arth–{M}oon trajectories in the restricted four-body problem}.
\newblock {\it \bibinfo{journal}{Advances in Space Research}\/},  {\it
  \bibinfo{volume}{45}\/}\bibinfo{issue}{(3)}, \bibinfo{pages}{398--409}.
\bibitem[{Astakhov \& Farrelly(2004)}]{astakhov2004capture}
\bibinfo{author}{Astakhov, S.~A.}, \& \bibinfo{author}{Farrelly, D.}
  (\bibinfo{year}{2004}).
\newblock \bibinfo{title}{Capture and escape in the elliptic restricted
  three-body problem}.
\newblock {\it \bibinfo{journal}{Monthly Notices of the Royal Astronomical
  Society}\/},  {\it \bibinfo{volume}{354}\/}\bibinfo{issue}{(4)},
  \bibinfo{pages}{971--979}.
\bibitem[{Bartsch et~al.(2008)Bartsch, Moix, Hernandez, Kawai \&
  Uzer}]{bartsch2008time}
\bibinfo{author}{Bartsch, T.}, \bibinfo{author}{Moix, J.~M.},
  \bibinfo{author}{Hernandez, R.}, \bibinfo{author}{Kawai, S.}, \&
  \bibinfo{author}{Uzer, T.} (\bibinfo{year}{2008}).
\newblock \bibinfo{title}{Time-dependent transition state theory}.
\newblock In \bibinfo{editor}{S.~A. Rice} (Ed.), {\it
  \bibinfo{booktitle}{Advances in Chemical Physics}\/} (p.
  \bibinfo{pages}{191}).
\newblock \bibinfo{publisher}{John Wiley \& Sons, Inc.} volume
  \bibinfo{volume}{140}.
\bibitem[{Bihan et~al.(2017)Bihan, Masdemont, G{\'{o}}mez \&
  Lizy-Destrez}]{Le_Bihan_2017}
\bibinfo{author}{Bihan, B.~L.}, \bibinfo{author}{Masdemont, J.~J.},
  \bibinfo{author}{G{\'{o}}mez, G.}, \& \bibinfo{author}{Lizy-Destrez, S.}
  (\bibinfo{year}{2017}).
\newblock \bibinfo{title}{Invariant manifolds of a non-autonomous
  quasi-bicircular problem computed via the parameterization method}.
\newblock {\it \bibinfo{journal}{Nonlinearity}\/},  {\it
  \bibinfo{volume}{30}\/}\bibinfo{issue}{(8)}, \bibinfo{pages}{3040--3075}.
\bibitem[{Broucke(1969)}]{broucke1969stability}
\bibinfo{author}{Broucke, R.} (\bibinfo{year}{1969}).
\newblock \bibinfo{title}{Stability of periodic orbits in the elliptic,
  restricted three-body problem}.
\newblock {\it \bibinfo{journal}{AIAA Journal}\/},  {\it
  \bibinfo{volume}{7}\/}\bibinfo{issue}{(6)}, \bibinfo{pages}{1003--1009}.
\bibitem[{Condon \& Pearson(2001)}]{CoPe2001}
\bibinfo{author}{Condon, G.~L.}, \& \bibinfo{author}{Pearson, D.~P.}
  (\bibinfo{year}{2001}).
\newblock \bibinfo{title}{The role of humans in libration point missions with
  specific application to an {E}arth-{M}oon libration point gateway station}.
\newblock In {\it \bibinfo{booktitle}{AAS/AIAA Astrodynamics Specialist
  Conference}\/}.
\newblock \bibinfo{address}{Quebec City, Canada}.
\newblock \bibinfo{note}{{P}aper {N}o. {AAS} 01-307}.
\bibitem[{Conley(1968)}]{Conley1968}
\bibinfo{author}{Conley, C.~C.} (\bibinfo{year}{1968}).
\newblock \bibinfo{title}{Low energy transit orbits in the restricted
  three-body problem}.
\newblock {\it \bibinfo{journal}{SIAM J. Appl. Math.}\/},  {\it
  \bibinfo{volume}{16}\/}, \bibinfo{pages}{732--746}.
\bibitem[{Conley(1969)}]{Conley1969}
\bibinfo{author}{Conley, C.~C.} (\bibinfo{year}{1969}).
\newblock \bibinfo{title}{On the ultimate behavior of orbits with respect to an
  unstable critical point. {I}. {O}scillating, asymptotic, and capture orbits}.
\newblock {\it \bibinfo{journal}{J. Differential Equations}\/},  {\it
  \bibinfo{volume}{5}\/}, \bibinfo{pages}{136--158}.
\bibitem[{Cronin et~al.(1964)Cronin, Richards \& Russell}]{cronin1964}
\bibinfo{author}{Cronin, J.}, \bibinfo{author}{Richards, P.}, \&
  \bibinfo{author}{Russell, L.} (\bibinfo{year}{1964}).
\newblock \bibinfo{title}{Some periodic solutions of a four-body problem}.
\newblock {\it \bibinfo{journal}{Icarus}\/},  {\it \bibinfo{volume}{3}\/},
  \bibinfo{pages}{423--428}.
\bibitem[{Dellnitz et~al.(2005)Dellnitz, Junge, Lo, Marsden, Padberg, Preis,
  Ross \& Thiere}]{DeJuLoMaPaPrRoTh2005}
\bibinfo{author}{Dellnitz, M.}, \bibinfo{author}{Junge, O.},
  \bibinfo{author}{Lo, M.~W.}, \bibinfo{author}{Marsden, J.~E.},
  \bibinfo{author}{Padberg, K.}, \bibinfo{author}{Preis, R.},
  \bibinfo{author}{Ross, S.~D.}, \& \bibinfo{author}{Thiere, B.}
  (\bibinfo{year}{2005}).
\newblock \bibinfo{title}{Transport of {M}ars-crossing asteroids from the
  quasi-{H}ilda region}.
\newblock {\it \bibinfo{journal}{Physical Review Letters}\/},  {\it
  \bibinfo{volume}{94}\/}, \bibinfo{pages}{231102}.
\bibitem[{Gawlik et~al.(2009)Gawlik, Marsden, {Du {T}oit} \&
  Campagnola}]{GaMaDuCa2009}
\bibinfo{author}{Gawlik, E.~S.}, \bibinfo{author}{Marsden, J.~E.},
  \bibinfo{author}{{Du {T}oit}, P.~C.}, \& \bibinfo{author}{Campagnola, S.}
  (\bibinfo{year}{2009}).
\newblock \bibinfo{title}{Lagrangian coherent structures in the planar elliptic
  restricted three-body problem}.
\newblock {\it \bibinfo{journal}{Celestial Mechanics and Dynamical
  Astronomy}\/},  {\it \bibinfo{volume}{103}\/}, \bibinfo{pages}{227--249}.
\bibitem[{G\'omez et~al.(2004)G\'omez, Koon, Lo, Marsden, Masdemont \&
  Ross}]{GoKoLoMaMaRo2004}
\bibinfo{author}{G\'omez, G.}, \bibinfo{author}{Koon, W.~S.},
  \bibinfo{author}{Lo, M.~W.}, \bibinfo{author}{Marsden, J.~E.},
  \bibinfo{author}{Masdemont, J.}, \& \bibinfo{author}{Ross, S.~D.}
  (\bibinfo{year}{2004}).
\newblock \bibinfo{title}{Connecting orbits and invariant manifolds in the
  spatial three-body problem}.
\newblock {\it \bibinfo{journal}{Nonlinearity}\/},  {\it
  \bibinfo{volume}{17}\/}, \bibinfo{pages}{1571--1606}.
\bibitem[{G{\'o}mez et~al.(2003)G{\'o}mez, Masdemont \&
  Mondelo}]{gomez2003dynamical}
\bibinfo{author}{G{\'o}mez, G.}, \bibinfo{author}{Masdemont, J.}, \&
  \bibinfo{author}{Mondelo, J.} (\bibinfo{year}{2003}).
\newblock \bibinfo{title}{Dynamical substitutes of the libration points for
  simplified solar system models}.
\newblock In {\it \bibinfo{booktitle}{Libration point orbits and
  applications}\/} (pp. \bibinfo{pages}{373--397}).
\newblock \bibinfo{publisher}{World Scientific}.
\bibitem[{Guckenheimer \& Holmes(2013)}]{guckenheimer2013nonlinear}
\bibinfo{author}{Guckenheimer, J.}, \& \bibinfo{author}{Holmes, P.}
  (\bibinfo{year}{2013}).
\newblock {\it \bibinfo{title}{Nonlinear oscillations, dynamical systems, and
  bifurcations of vector fields}\/} volume~\bibinfo{volume}{42}.
\newblock \bibinfo{publisher}{Springer Science \& Business Media}.
\bibitem[{Guo \& Lei(2019)}]{guo2019families}
\bibinfo{author}{Guo, Q.}, \& \bibinfo{author}{Lei, H.} (\bibinfo{year}{2019}).
\newblock \bibinfo{title}{Families of {E}arth--{M}oon trajectories with
  applications to transfers towards {S}un--{E}arth libration point orbits}.
\newblock {\it \bibinfo{journal}{Astrophysics and Space Science}\/},  {\it
  \bibinfo{volume}{364}\/}\bibinfo{issue}{(3)}, \bibinfo{pages}{1--9}.
\bibitem[{Jaff\'e et~al.(2002)Jaff\'e, Ross, Lo, Marsden, Farrelly \&
  Uzer}]{JaRoLoMaFaUz2002}
\bibinfo{author}{Jaff\'e, C.}, \bibinfo{author}{Ross, S.~D.},
  \bibinfo{author}{Lo, M.~W.}, \bibinfo{author}{Marsden, J.~E.},
  \bibinfo{author}{Farrelly, D.}, \& \bibinfo{author}{Uzer, T.}
  (\bibinfo{year}{2002}).
\newblock \bibinfo{title}{Theory of asteroid escape rates}.
\newblock {\it \bibinfo{journal}{Physical Review Letters}\/},  {\it
  \bibinfo{volume}{89}\/}, \bibinfo{pages}{011101}.
\bibitem[{Jorba et~al.(2020)Jorba, Jorba-Cusc\'o \& Rosales}]{jorbaetal}
\bibinfo{author}{Jorba, A.}, \bibinfo{author}{Jorba-Cusc\'o, M.}, \&
  \bibinfo{author}{Rosales, J.~J.} (\bibinfo{year}{2020}).
\newblock \bibinfo{title}{The vicinity of the {E}arth–{M}oon ${L}_1$ point in
  the bicircular problem}.
\newblock {\it \bibinfo{journal}{Celestial Mechanics and Dynamical
  Astronomy}\/},  {\it \bibinfo{volume}{132}\/}\bibinfo{issue}{(2)}.
\bibitem[{Jorba \& Nicol{\'a}s(2020)}]{jorba2020transport}
\bibinfo{author}{Jorba, A.}, \& \bibinfo{author}{Nicol{\'a}s, B.}
  (\bibinfo{year}{2020}).
\newblock \bibinfo{title}{Transport and invariant manifolds near {L}3 in the
  {E}arth-{M}oon {B}icircular model}.
\newblock {\it \bibinfo{journal}{Communications in Nonlinear Science and
  Numerical Simulation}\/},  {\it \bibinfo{volume}{89}\/},
  \bibinfo{pages}{105327}.
\bibitem[{Jordan \& Smith(2007)}]{jordan2007nonlinear}
\bibinfo{author}{Jordan, D.~W.}, \& \bibinfo{author}{Smith, P.}
  (\bibinfo{year}{2007}).
\newblock {\it \bibinfo{title}{Nonlinear ordinary differential equations: an
  introduction for scientists and engineers}\/} volume~\bibinfo{volume}{10}.
\newblock \bibinfo{publisher}{Oxford University Press}.
\bibitem[{Koon et~al.(2000)Koon, Lo, Marsden \& Ross}]{KoLoMaRo2000}
\bibinfo{author}{Koon, W.~S.}, \bibinfo{author}{Lo, M.~W.},
  \bibinfo{author}{Marsden, J.~E.}, \& \bibinfo{author}{Ross, S.~D.}
  (\bibinfo{year}{2000}).
\newblock \bibinfo{title}{Heteroclinic connections between periodic orbits and
  resonance transitions in celestial mechanics}.
\newblock {\it \bibinfo{journal}{Chaos}\/},  {\it \bibinfo{volume}{10}\/},
  \bibinfo{pages}{427--469}.
\bibitem[{Koon et~al.(2001{\natexlab{a}})Koon, Lo, Marsden \&
  Ross}]{KoLoMaRo2001a}
\bibinfo{author}{Koon, W.~S.}, \bibinfo{author}{Lo, M.~W.},
  \bibinfo{author}{Marsden, J.~E.}, \& \bibinfo{author}{Ross, S.~D.}
  (\bibinfo{year}{2001}{\natexlab{a}}).
\newblock \bibinfo{title}{Low energy transfer to the {M}oon}.
\newblock {\it \bibinfo{journal}{Celestial Mechanics and Dynamical
  Astronomy}\/},  {\it \bibinfo{volume}{81}\/}, \bibinfo{pages}{63--73}.
\bibitem[{Koon et~al.(2001{\natexlab{b}})Koon, Lo, Marsden \&
  Ross}]{KoLoMaRo2001}
\bibinfo{author}{Koon, W.~S.}, \bibinfo{author}{Lo, M.~W.},
  \bibinfo{author}{Marsden, J.~E.}, \& \bibinfo{author}{Ross, S.~D.}
  (\bibinfo{year}{2001}{\natexlab{b}}).
\newblock \bibinfo{title}{Resonance and capture of {J}upiter comets}.
\newblock {\it \bibinfo{journal}{Celestial Mechanics and Dynamical
  Astronomy}\/},  {\it \bibinfo{volume}{81}\/}, \bibinfo{pages}{27--38}.
\bibitem[{Koon et~al.(2011)Koon, Lo, Marsden \& Ross}]{kolomaro}
\bibinfo{author}{Koon, W.~S.}, \bibinfo{author}{Lo, M.~W.},
  \bibinfo{author}{Marsden, J.~E.}, \& \bibinfo{author}{Ross, S.~D.}
  (\bibinfo{year}{2011}).
\newblock {\it \bibinfo{title}{Dynamical Systems, the Three-Body Problem and
  Space Mission Design}\/}.
\newblock \bibinfo{publisher}{Marsden Books, ISBN 978-0-615-24095-4}.
\bibitem[{Kraj{\v{n}}{\'a}k \& Waalkens(2018)}]{krajvnak2018phase}
\bibinfo{author}{Kraj{\v{n}}{\'a}k, V.}, \& \bibinfo{author}{Waalkens, H.}
  (\bibinfo{year}{2018}).
\newblock \bibinfo{title}{The phase space geometry underlying roaming reaction
  dynamics}.
\newblock {\it \bibinfo{journal}{Journal of Mathematical Chemistry}\/},  {\it
  \bibinfo{volume}{56}\/}\bibinfo{issue}{(8)}, \bibinfo{pages}{2341--2378}.
\bibitem[{Kumar et~al.(2021)Kumar, Anderson \& de~la Llave}]{kumar2021using}
\bibinfo{author}{Kumar, B.}, \bibinfo{author}{Anderson, R.~L.}, \&
  \bibinfo{author}{de~la Llave, R.} (\bibinfo{year}{2021}).
\newblock \bibinfo{title}{Using gpus and the parameterization method for rapid
  search and refinement of connections between tori in periodically perturbed
  planar circular restricted 3-body problems}.
\newblock {\it \bibinfo{journal}{arXiv preprint arXiv:2109.14814}\/}, .
\bibitem[{Lichtenberg \& Lieberman(1992)}]{lichtenberg1992regular}
\bibinfo{author}{Lichtenberg, A.~J.}, \& \bibinfo{author}{Lieberman, M.~A.}
  (\bibinfo{year}{1992}).
\newblock {\it \bibinfo{title}{Regular and Chaotic Dynamics}\/}
  volume~\bibinfo{volume}{38}.
\newblock (\bibinfo{edition}{2nd} ed.).
\newblock \bibinfo{publisher}{Springer}.
\bibitem[{Llibre et~al.(1985)Llibre, Martinez \& Sim\'o}]{LlMaSi1985}
\bibinfo{author}{Llibre, J.}, \bibinfo{author}{Martinez, R.}, \&
  \bibinfo{author}{Sim\'o, C.} (\bibinfo{year}{1985}).
\newblock \bibinfo{title}{Transversality of the invariant manifolds associated
  to the {L}yapunov family of periodic orbits near {L}2 in the restricted
  three-body problem}.
\newblock {\it \bibinfo{journal}{J. Diff. Eqns.}\/},  {\it
  \bibinfo{volume}{58}\/}, \bibinfo{pages}{104--156}.
\bibitem[{Lo \& Ross(2001)}]{lo2001lunar}
\bibinfo{author}{Lo, M.~W.}, \& \bibinfo{author}{Ross, S.~D.}
  (\bibinfo{year}{2001}).
\newblock \bibinfo{title}{The lunar {L}1 gateway: portal to the stars and
  beyond}.
\newblock In {\it \bibinfo{booktitle}{AIAA Space 2001 Conference and
  Exposition}\/} (p. \bibinfo{pages}{4768}).
\bibitem[{MacKay(1990)}]{MacKay1990}
\bibinfo{author}{MacKay, R.~S.} (\bibinfo{year}{1990}).
\newblock \bibinfo{title}{Flux over a saddle}.
\newblock {\it \bibinfo{journal}{Physics Letters A}\/},  {\it
  \bibinfo{volume}{145}\/}, \bibinfo{pages}{425--427}.
\bibitem[{Marsden \& Ratiu(1999)}]{MaRa1999}
\bibinfo{author}{Marsden, J.~E.}, \& \bibinfo{author}{Ratiu, T.~S.}
  (\bibinfo{year}{1999}).
\newblock {\it \bibinfo{title}{Introduction to Mechanics and Symmetry}\/}
  volume~\bibinfo{volume}{17} of {\it \bibinfo{series}{Texts in Applied
  Mathematics}\/}.
\newblock \bibinfo{address}{New York}: \bibinfo{publisher}{Springer-Verlag}.
\bibitem[{McCarthy \& Howell(2020)}]{mccarthy2020cislunar}
\bibinfo{author}{McCarthy, B.~P.}, \& \bibinfo{author}{Howell, K.~C.}
  (\bibinfo{year}{2020}).
\newblock \bibinfo{title}{Cislunar transfer design exploiting periodic and
  quasi-periodic orbital structures in the four-body problem}.
\newblock In {\it \bibinfo{booktitle}{71st International Astronautical
  Congress, Virtual}\/}.
\bibitem[{McGehee(1969{\natexlab{a}})}]{McGehee1969}
\bibinfo{author}{McGehee, R.} (\bibinfo{year}{1969}{\natexlab{a}}).
\newblock {\it \bibinfo{title}{Some homoclinic orbits for the restricted
  three-body problem}\/}.
\newblock Ph.D. thesis University of Wisconsin, Madison.
\bibitem[{McGehee(1969{\natexlab{b}})}]{mcgehee1969some}
\bibinfo{author}{McGehee, R.~P.} (\bibinfo{year}{1969}{\natexlab{b}}).
\newblock {\it \bibinfo{title}{Some homoclinic orbits for the restricted
  three-body problem}\/}.
\newblock \bibinfo{publisher}{The University of Wisconsin-Madison}.
\bibitem[{Moser(1958)}]{Moser1958}
\bibinfo{author}{Moser, J.} (\bibinfo{year}{1958}).
\newblock \bibinfo{title}{On the generalization of a theorem of {L}iapunov}.
\newblock {\it \bibinfo{journal}{Comm. Pure Appl. Math.}\/},  {\it
  \bibinfo{volume}{11}\/}, \bibinfo{pages}{257--271}.
\bibitem[{Moser(1973)}]{Moser1973}
\bibinfo{author}{Moser, J.} (\bibinfo{year}{1973}).
\newblock {\it \bibinfo{title}{Stable and Random Motions in Dynamical Systems
  with Special Emphasis on Celestial Mechanics}\/}.
\newblock \bibinfo{publisher}{Princeton University Press}.
\bibitem[{Naik \& Ross(2017)}]{NaikRoss2017a}
\bibinfo{author}{Naik, S.}, \& \bibinfo{author}{Ross, S.~D.}
  (\bibinfo{year}{2017}).
\newblock \bibinfo{title}{Geometry of escaping dynamics in nonlinear ship
  motion}.
\newblock {\it \bibinfo{journal}{Communications in Nonlinear Science and
  Numerical Simulation}\/},  {\it \bibinfo{volume}{47}\/}, \bibinfo{pages}{48
  -- 70}.
\bibitem[{de~Oliveira et~al.(2020)de~Oliveira, Sousa-Silva \&
  Caldas}]{deoliveira2020order}
\bibinfo{author}{de~Oliveira, V.~M.}, \bibinfo{author}{Sousa-Silva, P.~A.}, \&
  \bibinfo{author}{Caldas, I.~L.} (\bibinfo{year}{2020}).
\newblock \bibinfo{title}{Order-chaos-order and invariant manifolds in the
  bounded planar {E}arth--{M}oon system}.
\newblock {\it \bibinfo{journal}{Celestial Mechanics and Dynamical
  Astronomy}\/},  {\it \bibinfo{volume}{132}\/}\bibinfo{issue}{(11)},
  \bibinfo{pages}{1--17}.
\bibitem[{Onozaki et~al.(2017)Onozaki, Yoshimura \& Ross}]{onozaki2017tube}
\bibinfo{author}{Onozaki, K.}, \bibinfo{author}{Yoshimura, H.}, \&
  \bibinfo{author}{Ross, S.~D.} (\bibinfo{year}{2017}).
\newblock \bibinfo{title}{Tube dynamics and low energy {E}arth-{M}oon transfers
  in the 4-body system}.
\newblock {\it \bibinfo{journal}{Advances in Space Research}\/},  {\it
  \bibinfo{volume}{60}\/}\bibinfo{issue}{(10)}, \bibinfo{pages}{2117--2132}.
\bibitem[{Oshima et~al.(2017)Oshima, Topputo, Campagnola \&
  Yanao}]{oshima2017analysis}
\bibinfo{author}{Oshima, K.}, \bibinfo{author}{Topputo, F.},
  \bibinfo{author}{Campagnola, S.}, \& \bibinfo{author}{Yanao, T.}
  (\bibinfo{year}{2017}).
\newblock \bibinfo{title}{Analysis of medium-energy transfers to the {M}oon}.
\newblock {\it \bibinfo{journal}{Celestial Mechanics and Dynamical
  Astronomy}\/},  {\it \bibinfo{volume}{127}\/}\bibinfo{issue}{(3)},
  \bibinfo{pages}{285--300}.
\bibitem[{Oshima \& Yanao(2014)}]{oshima2014applications}
\bibinfo{author}{Oshima, K.}, \& \bibinfo{author}{Yanao, T.}
  (\bibinfo{year}{2014}).
\newblock \bibinfo{title}{Applications of gravity assists in the bi-circular
  and bielliptic restricted four-body problem}.
\newblock In {\it \bibinfo{booktitle}{AAS/AIAA Space Flight Mechanics
  Meeting}\/}.
\newblock \bibinfo{address}{Santa Fe, New Mexico}.
\newblock \bibinfo{note}{{P}aper {N}o. {AAS} 14–234}.
\bibitem[{Paez \& Guzzo(2021)}]{Paez_2021}
\bibinfo{author}{Paez, R.~I.}, \& \bibinfo{author}{Guzzo, M.}
  (\bibinfo{year}{2021}).
\newblock \bibinfo{title}{Transits close to the lagrangian solutions $l_1$,
  $l_2$ in the elliptic restricted three-body problem}.
\newblock {\it \bibinfo{journal}{Nonlinearity}\/},  {\it
  \bibinfo{volume}{34}\/}\bibinfo{issue}{(9)}, \bibinfo{pages}{6417--6449}.
  \DOIprefix\doi{10.1088/1361-6544/ac13be}.
\bibitem[{Parker \& Chua(1989)}]{PaCh1989}
\bibinfo{author}{Parker, T.~S.}, \& \bibinfo{author}{Chua, L.~O.}
  (\bibinfo{year}{1989}).
\newblock {\it \bibinfo{title}{Practical Numerical Algorithms for Chaotic
  Systems}\/}.
\newblock \bibinfo{address}{New York}: \bibinfo{publisher}{Springer-Verlag}.
\bibitem[{Reddy(2008)}]{Reddy2008}
\bibinfo{author}{Reddy, F.} (\bibinfo{year}{2008}).
\newblock \bibinfo{title}{How scientists discovered a solar system
  ``superhighway''}.
\newblock {\it \bibinfo{journal}{Astronomy}\/},  {\it \bibinfo{volume}{36}\/},
  \bibinfo{pages}{38--43}.
\bibitem[{Ren \& Shan(2012)}]{ren2012numerical}
\bibinfo{author}{Ren, Y.}, \& \bibinfo{author}{Shan, J.}
  (\bibinfo{year}{2012}).
\newblock \bibinfo{title}{Numerical study of the three-dimensional transit
  orbits in the circular restricted three-body problem}.
\newblock {\it \bibinfo{journal}{Celestial Mechanics and Dynamical
  Astronomy}\/},  {\it \bibinfo{volume}{114}\/}\bibinfo{issue}{(4)},
  \bibinfo{pages}{415--428}.
\bibitem[{Ross(2006)}]{Ross2006}
\bibinfo{author}{Ross, S.~D.} (\bibinfo{year}{2006}).
\newblock \bibinfo{title}{The interplanetary transport network}.
\newblock {\it \bibinfo{journal}{American Scientist}\/},  {\it
  \bibinfo{volume}{94}\/}, \bibinfo{pages}{230--237}.
\bibitem[{Ross \& Scheeres(2007)}]{RoSc2007}
\bibinfo{author}{Ross, S.~D.}, \& \bibinfo{author}{Scheeres, D.~J.}
  (\bibinfo{year}{2007}).
\newblock \bibinfo{title}{Multiple gravity assists, capture, and escape in the
  restricted three-body problem}.
\newblock {\it \bibinfo{journal}{SIAM Journal on Applied Dynamical Systems}\/},
   {\it \bibinfo{volume}{6}\/}\bibinfo{issue}{(3)}, \bibinfo{pages}{576--596}.
  \DOIprefix\doi{10.1137/060663374}.
\bibitem[{Sim\'o et~al.(1995)Sim\'o, G\'omez, Jorba \&
  Masdemont}]{SiGoJoMa1995}
\bibinfo{author}{Sim\'o, C.}, \bibinfo{author}{G\'omez, G.},
  \bibinfo{author}{Jorba, A.}, \& \bibinfo{author}{Masdemont, J.}
  (\bibinfo{year}{1995}).
\newblock \bibinfo{title}{The bicircular model near the triangular libration
  points}.
\newblock In \bibinfo{editor}{A.~E. Roy}, \& \bibinfo{editor}{B.~A. Steves}
  (Eds.), {\it \bibinfo{booktitle}{From Newton to Chaos}\/} (pp.
  \bibinfo{pages}{343--370}).
\newblock \bibinfo{address}{New York}: \bibinfo{publisher}{Plenum Press}.
\bibitem[{Szebehely(1967)}]{Szebehely1967}
\bibinfo{author}{Szebehely, V.} (\bibinfo{year}{1967}).
\newblock {\it \bibinfo{title}{Theory of Orbits: The Restricted Problem of
  Three Bodies}\/}.
\newblock \bibinfo{address}{New York}: \bibinfo{publisher}{Academic}.
\bibitem[{Todorovi{\'c} et~al.(2020)Todorovi{\'c}, Wu \&
  Rosengren}]{todorovic2020arches}
\bibinfo{author}{Todorovi{\'c}, N.}, \bibinfo{author}{Wu, D.}, \&
  \bibinfo{author}{Rosengren, A.~J.} (\bibinfo{year}{2020}).
\newblock \bibinfo{title}{The arches of chaos in the solar system}.
\newblock {\it \bibinfo{journal}{Science Advances}\/},  {\it
  \bibinfo{volume}{6}\/}\bibinfo{issue}{(48)}, \bibinfo{pages}{eabd1313}.
\bibitem[{Topputo(2013)}]{topputo2013optimal}
\bibinfo{author}{Topputo, F.} (\bibinfo{year}{2013}).
\newblock \bibinfo{title}{On optimal two-impulse {E}arth--{M}oon transfers in a
  four-body model}.
\newblock {\it \bibinfo{journal}{Celestial Mechanics and Dynamical
  Astronomy}\/},  {\it \bibinfo{volume}{117}\/}\bibinfo{issue}{(3)},
  \bibinfo{pages}{279--313}.
\bibitem[{Waalkens \& Wiggins(2004)}]{waalkens2004direct}
\bibinfo{author}{Waalkens, H.}, \& \bibinfo{author}{Wiggins, S.}
  (\bibinfo{year}{2004}).
\newblock \bibinfo{title}{Direct construction of a dividing surface of minimal
  flux for multi-degree-of-freedom systems that cannot be recrossed}.
\newblock {\it \bibinfo{journal}{Journal of Physics A: Mathematical and
  General}\/},  {\it \bibinfo{volume}{37}\/}\bibinfo{issue}{(35)},
  \bibinfo{pages}{L435}.
\bibitem[{Wiggins(2003)}]{Wiggins2003}
\bibinfo{author}{Wiggins, S.} (\bibinfo{year}{2003}).
\newblock {\it \bibinfo{title}{Introduction to Applied Nonlinear Dynamical
  Systems and Chaos}\/} volume~\bibinfo{volume}{2} of {\it
  \bibinfo{series}{Texts in Applied Mathematics Science}\/}.
\newblock (\bibinfo{edition}{2nd} ed.).
\newblock \bibinfo{address}{Berlin}: \bibinfo{publisher}{Springer-Verlag}.
\bibitem[{Wu \& McCue(2008)}]{wu2008application}
\bibinfo{author}{Wu, W.}, \& \bibinfo{author}{McCue, L.}
  (\bibinfo{year}{2008}).
\newblock \bibinfo{title}{Application of the extended {M}elnikov's method for
  single-degree-of-freedom vessel roll motion}.
\newblock {\it \bibinfo{journal}{Ocean Engineering}\/},  {\it
  \bibinfo{volume}{35}\/}\bibinfo{issue}{(17-18)}, \bibinfo{pages}{1739--1746}.
\bibitem[{Zhong \& Ross(2020)}]{zhong2020geometry}
\bibinfo{author}{Zhong, J.}, \& \bibinfo{author}{Ross, S.~D.}
  (\bibinfo{year}{2020}).
\newblock \bibinfo{title}{Geometry of escape and transition dynamics in the
  presence of dissipative and gyroscopic forces in two degree of freedom
  systems}.
\newblock {\it \bibinfo{journal}{Communications in Nonlinear Science and
  Numerical Simulation}\/},  {\it \bibinfo{volume}{82}\/},
  \bibinfo{pages}{105033}.
\bibitem[{Zhong \& Ross(2021)}]{zhong2021global}
\bibinfo{author}{Zhong, J.}, \& \bibinfo{author}{Ross, S.~D.}
  (\bibinfo{year}{2021}).
\newblock \bibinfo{title}{Global invariant manifolds delineating transition and
  escape dynamics in dissipative systems: an application to snap-through
  buckling}.
\newblock {\it \bibinfo{journal}{Nonlinear Dynamics}\/},  {\it
  \bibinfo{volume}{104}\/}, \bibinfo{pages}{3109--3137}.

\end{thebibliography}

\appendix
\section{Initial Conditions}

In the BCP as described in Section \ref{BCP}, the Lagrange periodic orbit replacing the Earth-Moon \(L_1\) point has initial condition,
\[
\mathbf{\bar x}=
\begin{bmatrix}
   \bar x \\
   \bar y \\
   \bar p_x \\
   \bar p_y
\end{bmatrix}
=
\begin{bmatrix}
   0.837595408485656 \\
   0 \\
   0 \\
   0.827678389393936
\end{bmatrix}
\]
in the four-dimensional position-momentum phase space at phase $\theta = 0$.


In the ER3BP as described in Section \ref{ER3BP}, the Lagrange periodic orbit replacing the Earth-Moon \(L_1\) point has initial condition,
\[
\mathbf{\bar x}=
\begin{bmatrix}
   \bar x \\
   \bar y \\
   \bar p_x \\
   \bar p_y
\end{bmatrix}
=
\begin{bmatrix}
   0.792718947200736 \\
   0 \\
   0.000001145970495 \\
   0.886145419995798
\end{bmatrix}
\]
in the four-dimensional position-momentum phase space at phase $\theta = 0$. We suspect that $\bar p_x \ne 0$ is a numerical artifact.

\section{Proof of Proposition 1}

\begin{proof}
The assumed quadratic Hamiltonian function is,
\begin{equation}
    \Tilde{H}_2(\mathbf{x}) = \Tilde{H}_2(q_1,p_1,q_2,p_2) = \Tilde{\lambda} q_1 p_1 + \tfrac{1}{2} \Tilde{\nu}( q_2^2 + p_2^2 ).
\end{equation}
Hamilton's canonical equations generated by this Hamiltonian are linear,
\begin{equation}
\begin{split}
\dot{\mathbf{x}} 
=J \nabla \Tilde{H}_2(\mathbf{x})
 &=
\begin{bmatrix}
 0 &  1 & 0 &  0 \\
-1 &  0 & 0 &  0 \\
 0 &  0 & 0 &  1 \\
 0 &  0 & -1 & 0
\end{bmatrix}
\begin{bmatrix}
    \Tilde{\lambda} p_1 \\
    \Tilde{\lambda} q_1 \\
    \Tilde{\nu}  q_2 \\
    \Tilde{\nu}  p_2 \\
\end{bmatrix} \\
&=
\underbrace{\begin{bmatrix}
 \Tilde{\lambda} &  0                & 0                & 0              \\
 0               &  -\Tilde{\lambda} & 0                & 0              \\
 0               &  0                & 0                & \Tilde{\nu} \\
 0               &  0                & -\Tilde{\nu}  & 0
\end{bmatrix}}_{\mathbf{A}}
\mathbf{x}.
\end{split}
\label{eq:linearquadsys}
\end{equation}
which is of the form (\ref{linearized_pcr3bp_eigen}) with $\lambda=\tilde{\lambda}$, $\nu=\tilde{\nu}$, where $\mathbf{x}=(q_1,p_1,q_2,p_2)^T$.
\\\\
\indent It is straightforward to show analytically that the solution to the linear differential equation \eqref{eq:linearquadsys} is, 
\begin{equation}
\begin{split}
    &\mathbf{x}(t) = 
    e^{\mathbf{A}t} \mathbf{x}(0)=
    \begin{bmatrix}
    e^{ \Tilde{\lambda} t} &   0           &   ~~0  & 0 \\
    0      &   e^{ -\Tilde{\lambda} t} &   ~~0  & 0 \\
    0      &   0           &   ~~\cos{(\Tilde{\nu} t)}  & \sin{(\Tilde{\nu} t)} \\
    0      &   0           &  -\sin{(\Tilde{\nu} t)}  & \cos{(\Tilde{\nu} t)}
\end{bmatrix}
 \mathbf{x}(0), \\
 &\quad {\rm where}
\quad
    \mathbf{x}(0) = 
\begin{bmatrix}
    q_{1_0} \\
    p_{1_0} \\
    q_{2_0} \\
    p_{2_0}
\end{bmatrix}.
\label{eq:linquadsysic}
\end{split}
\end{equation}
We note that $e^{\mathbf{A}T}$ is of the form $\mathbf{\Lambda}$ from (\ref{eq:lambda_form}) with 
\begin{equation}
    \sigma = e^{ \Tilde{\lambda} T}, \quad \psi = \Tilde{\nu}T, 
\end{equation}
which is equivalent to  (\ref{effective_params}).
Therefore,
\begin{equation}
\mathbf{x}(T)=\mathbf{\Lambda}\mathbf{x}(0)
\end{equation}
And thus \(\Tilde{H}_2(\mathbf{x})\) generates the linear symplectic map \(\mathbf{x} \mapsto \mathbf{\Lambda} \mathbf{x}\), with $ \mathbf{\Lambda}$ as in \eqref{eq:lambda_form}. \qed
\end{proof}

\section{Continuation Visualization for the ER3BP $L_1$ Lagrange Periodic Orbit}\label{appendix:c}

\begin{figure}[]
\begin{center}
\includegraphics[width=0.6\textwidth]{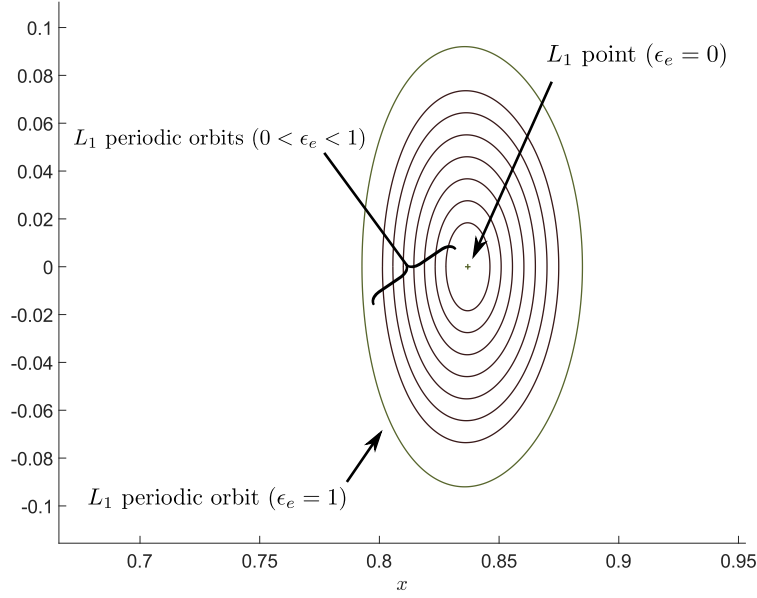}
\caption{A family of periodic orbits for different eccentricities emanating from the unperturbed $L_1$ point in the ER3BP.}
\label{fig:ER3BPcontinuation}
\end{center}
\vspace{-8mm}
\end{figure}

The ER3BP $L_1$ Lagrange periodic orbit can be obtained through continuation using the methodology described in Section 3.4. Let the true eccentricity of the system be $e$. A rescaled eccentricity is given by $e \epsilon_e$ where $\epsilon_e = 0$ for 0 eccentricity and $\epsilon_e = 1$ for the true eccentricity. Substituting $e \epsilon_e$ into the equations of motion and slowly increasing $\epsilon_e$ while refining the Lagrange periodic orbit for each perturbation of the parameter demonstrates continuity between the Lagrange point and the full eccentricity Lagrange periodic orbit (see Figure \ref{fig:ER3BPcontinuation}).

\end{document}